\documentclass{article}

\usepackage{amsthm}
\usepackage{amsbsy}
\usepackage{amscd}
\usepackage[dvips]{graphicx}
\usepackage{amsmath,amscd,amssymb,epsfig,subfig}
\usepackage[all]{xy}
\usepackage{verbatim}
\usepackage{graphics,graphicx}
\usepackage{tikz}
\newcommand{\fun}{\mathsf{F}}
\DeclareMathOperator{\End}{\mathrm{End}}

\DeclareMathOperator{\Ad}{Ad}

\newcommand{\cla}{c_b}
\newcommand{\COMPLEXS}{\mathbb C}

\newcommand{\EULER}{N}

\newcommand{\FUNCTOR}{\mathsf F}

\newcommand{\IMUN}{\mathsf i}
\newcommand{\imun}{\mathsf i}

\newcommand{\la}{b}

\newcommand{\MCG}{\mathcal M_\Sigma}
\newcommand{\MOM}{\mathsf p}

\newcommand{\NTR}{N}

\newcommand{\PERMUTE}{P}
\newcommand{\PGROUP}{S}
\newcommand{\POSI}{\mathsf q}
\newcommand{\PTOLEMY}{T}

\newcommand{\QDILOG}{\varphi_{b}}

\newcommand{\REALS}{\mathbb R}

\newcommand{\ROTATE}{\mathsf A}

\newcommand{\SDIT}{\Delta_\Sigma}

\newcommand{\sfa}{\mathsf a}

\newcommand{\sfu}{\mathsf u}

\newcommand{\SURFACE}{\Sigma}

\newcommand{\TRIANGLES}{T}

\newcommand{\vsp}{\mathcal H}

\newcommand{\gpd}{\mathcal{G}_\Sigma}

\theoremstyle{plain}
\newtheorem{theorem}{Theorem}[section]
\newtheorem{lemma}{Lemma}[section]
\newtheorem{corollary}[theorem]{Corollary}
\newtheorem{proposition}{Proposition}[section]

\newtheorem{definition}{Definition}[section]
\theoremstyle{remark}
\newtheorem{remark}{Remark}[section]

\newcommand{\Z}{{\mathbb{Z}}}
 
\newcommand{\C}{{\mathbb{C}}}

\begin{document}

\date{\today}

\title{Centrally extended mapping class groups from
quantum Teichm\"uller theory\footnote{
This version: {February 2010}.  L.F. was partially supported by
the ANR-06-BLAN-0311 Repsurf and ANR 2011 BS 01 020 01 ModGroup. 
R.M.K. is partially supported by Swiss National Science Foundation.}}
  \author{
\begin{tabular}{cc}
 Louis Funar &  Rinat M. Kashaev\\
\small \em Institut Fourier BP 74, UMR 5582
&\small \em Section de Math\'ematiques, Case Postale 64\\
\small \em University of Grenoble I &\small \em Universit\'e  de Gen\`eve\\
\small \em 38402 Saint-Martin-d'H\`eres cedex, France
&\small \em 2-4, rue du Li\`evre, 1211 Geneve 4, Suisse\\
\small \em e-mail: {\tt funar@fourier.ujf-grenoble.fr}
& \small \em e-mail: {\tt rinat.kashaev@unige.ch} \\
\end{tabular}
}

\maketitle

\begin{abstract}
The central extension of the mapping class groups of 
punctured surfaces of finite type that arises in
quantum Teichm\"uller theory is 12 times the Meyer class plus 
the Euler classes of the punctures. This is analogous
to the result obtained in \cite{FS}
for the Thompson groups.

\vspace{0.1cm}


\end{abstract}

\section*{Introduction}

The quantum theory of Teichm\"uller spaces of punctured surfaces of finite type, originally constructed in \cite{CF,Ka} and subsequently generalized to higher rank Lie groups and cluster algebras in  \cite{FG1,FG}, leads to one parameter families
of  projective unitary representations of Ptolemy modular groupoids associated
to ideal triangulations of punctured surfaces. We will call such representations (quantum) dilogarithmic representations, since the main ingredient in the theory is Faddeev's quantum dilogarithm function first introduced in the context of quantum integrable systems by L.D. Faddeev in \cite{Fad}.

\vspace{0.1cm}\noindent 
These representations are infinite
dimensional so that a priori it is not clear if they come from suitable $2+1$-dimensional topological quantum field theories (TQFT). Nonetheless, it is expected that in the singular limit, when the deformation parameter  tends to a root of unity\footnote{One should distinguish between two different limits, depending on whether $\frac{\log(q)}{2\pi i}$ tends to a positive or a negative rational number. In the case when this limit is a positive rational number, the limit of the representation is non-singular and so it stays infinite dimensional.}, the "renormalized" theory corresponds to a finite dimensional TQFT first defined in \cite{Ka3,Ka4} by using the cyclic representations of the Borel Hopf sub-algebra $BU_q(sl(2))$, and subsequently developed and generalized in \cite{BB}. One can get the same finite dimensional representations of Ptolemy modular groupoids directly from compact representations of quantum Teichm\"uller theory at roots of unity \cite{BL,BBL,Ka}.

\vspace{0.1cm}\noindent
Projective representations of a group are well known to be equivalent to representations of
central extensions of the same group by  means of the following procedure.
To a group $G$, a $\mathbb{C}$-vector space $V$ and a group homomorphism
$h:G\to PGL(V)\simeq GL(V)/\mathbb{C}^*$, where  $\mathbb{C}^*$ is identified with a (normal) subgroup of $GL(V)$ through the embedding $z\mapsto z\operatorname{id}_V$, one can associate a central extension $\widetilde{G}$
of $G$ by a sub-group $A$ of $\C^*$ together with a representation $\widetilde{h}:\widetilde{G}\to GL(V)$ such that the following diagram is commutative and has exact rows:
\[
\xymatrix{
1\ar[r]& \mathbb{C^*}\ar[r]& GL(V)\ar[r]&PGL(V)\ar[r] &1\\
1\ar[r]& A\ar[r]\ar@{^{(}->}[u]& \widetilde{G}\ar[r]\ar[u]_{\widetilde{h}}&G\ar[r]\ar[u]_h &1
}
\]

\vspace{0.1cm}\noindent 
One such extension is the pull-back $\mathbf{\widetilde{G}}$ of the 
central extension $GL(V)\to PGL(V)$ under the homomorphism $G\to PGL(V)$, which is canonically defined.   
However it is possible to find also smaller extensions associated to 
proper sub-groups $A\subset \C^*$.  
The central extension $\widetilde{G}$ associated to the 
smallest possible sub-group $A\subset \C^*$ 
for which there exists a linear representation as in the diagram above resolving the projective representation of $G$ will be called the {\em 
minimal reduction} of $\mathbf{\widetilde{G}}$.

\vspace{0.1cm}\noindent
In this light, quantum Teichm\"uller theory gives rise to representations of certain central extensions of the surface mapping class groups which are the vertex groups of the Ptolemy modular groupoids.
The main goal of this paper is to identify the isomorphism classes of those central extensions.
Namely, by using the quantization approach of \cite{Ka}, we extend the analysis of the particular case of a once punctured genus three surface performed in \cite{Ka2} to arbitrary punctured surfaces of finite type.

\vspace{0.1cm}\noindent
 Let a group $G$ with a given presentation be identified as the quotient group $F/R$, where $F$ is a free group and $R$, the normal
subgroup  generated by the relations. Then, a central extension of $G$ can be obtained from
a homomorphism $\overline{h}\colon F\to GL(V)$ with
the property $\overline{h}(R)\subset \C^*$ so that it induces a homomorphism $h\colon G\to PGL(V)$.
In this case,  the homomorphism $\overline{h}$ will be called
an {\em almost linear representation} of $G$, in order to distinguish
it from a  projective representation.

\vspace{0.1cm}\noindent
In quantum Teichm\"uller theory, central extensions of surface mapping class groups appear through almost linear representations. Specifically,
let $\Gamma^s_{g,r}$ be the mapping class group of a surface
$\Sigma^s_{g,r}$ of genus $g$ with $r$ boundary components and $s$ punctures.
These are mapping classes of homeomorphisms which
fix the boundary point-wise and fix the set of punctures (not necessarily point-wise).
Denoting $\Gamma^s_{g}=\Gamma^s_{g,0}$,
the projective representations of
$\Gamma^s_g$ for $(2g-2+2s)s >0$, constructed in \cite{Ka,Ka2}, are almost linear
representations corresponding to certain central extensions $\widetilde{\Gamma^s_{g}}$.
By considering embeddings $\Sigma^s_{g,r}\subset  \Sigma^{t}_{h,0}$,
the central extensions $\widetilde{\Gamma^s_{g}}$ can be used to define central extensions for the mapping class groups $\Gamma^s_{g,r}$ for $s\geq 1$. 
According to \cite{PR1}, any embedding $\Sigma^s_{g,r}\subset  \Sigma^{t}_{h,0}$,
for which $\Sigma^{t}_{h}\setminus \Sigma^s_{g,r}$
contains no disk, punctured disk or cylinder components, 
induces an embedding of the
corresponding mapping class groups. Using this fact,
we can define the central extension $\widetilde{\Gamma^s_{g,r}}$ as
the pull-back of the
central extension $\widetilde{\Gamma^t_{h}}$ by the injective
homomorphism
$\Gamma^s_{g,r}\hookrightarrow \Gamma^t_{h}$
induced by an embedding of the corresponding surfaces. A priori, it is not clear  whether such definition depends on a particular choice of the embedding, but our main result below shows
that this is indeed the case.

\vspace{0.1cm}\noindent
Central extensions by an Abelian group $A$ of a given group $G$ 
are known to be classified, up
to isomorphism,  by elements of the 2-cohomology group $H^2(G;A)$. 
In the case of surface mapping class groups $\Gamma_{g,r}^s$, the latter  was first computed by Harer in \cite{Harer} for $g\geq 5$
and further completed by Korkmaz and Stipsicz in \cite{KS} for $g\geq 4$ (see also \cite{Ko} for a survey).
Specifically, we have
\[
H^2(\Gamma^s_{g,r})=\Z^{s+1}, \: \mbox{\rm if }\, g\geq 4,
\]
where the generators are given by (one fourth of) the  Meyer signature
class $\chi$ (it is the only generator for the groups $H^2(\Gamma_g)\cong H^2(\Gamma_{g,1})\simeq\mathbb{Z}$,
see \cite{Me,Harer,KS} for its definition) and $s$ Euler classes $e_i$ associated with $s$ punctures. In the case when $g=3$, the group
$H^2(\Gamma^s_{3,r})$ still contains the sub-group $\Z^{s+1}$ generated
by the above mentioned classes, but it is not known whether
there are other (2-torsion) classes. When $g=2$ we will show that 
$H^2(\Gamma^s_{2,r})$ contains  the subgroup 
$\Z/10\Z\oplus\Z^s$, whose torsion 
part is generated by $\chi$ and whose free part is generated by the Euler 
classes.  The Universal Coefficients Theorem permits then to 
compute $H^2(G;A)$ for every Abelian group $A$.  
 
\vspace{0.2cm}\noindent 
Denote as above by $\mathbf{\widetilde{\Gamma^s_{g,r}}}$ the 
canonical central extension 
of $\Gamma_{g,r}^s$ by $\C^*$ which is obtained as the pull-back of 
the canonical central extension $GL(\mathcal H)\to PGL(\mathcal H)$ 
under the quantum projective representation 
 associated to a semi-symmetric $T$ 
in the Hilbert space $\mathcal H$ (see the next section). 
Quantum representations depend on some parameter $\zeta\in \C^*$. 
Our main result is the following theorem.

\begin{theorem}\label{dilog}
The central extension $\mathbf{\widetilde{\Gamma^s_{g,r}}}$ can be reduced 
to a minimal central extension $\widetilde{\Gamma^s_{g,r}}$ of 
${\Gamma}^s_{g,r}$ by  a cyclic Abelian $A\subset \C^*$, where $A$ is 
the subgroup of $\C^*$ generated by $\zeta^{-6}$. 
Moreover its  cohomology class is 
\[c_{{\widetilde{\Gamma^s_{g,r}} }}= 12 \chi +\sum_{i=1}^s e_i \in 
H^2(\Gamma^s_{g,r};A)\] 
if $g\geq 2$ and $s\geq 4$. 
Here $\chi$ and $e_i$ are one fourth of the Meyer signature class and 
respectively the $i$-th Euler class with $A$ coefficients.
\end{theorem}

\vspace{0.2cm}\noindent 
There is a  geometric interpretation of this extension. 
\begin{corollary}\label{cor1}
Let us consider the extension $\widehat{\Gamma_{g,r+s}}$ of class 
$12\chi$. Then there is an exact sequence: 
\[ 1\to A^{s-1}\to \widehat{\Gamma_{g,r+s}}\to \widetilde{\Gamma_{g,r}^s}\to 1\] 
In some sense the quantum representations of punctured mapping class groups 
can be lifted to the mapping class groups of surfaces with boundary 
obtained by blowing up the punctures. 
\end{corollary}

\begin{corollary}\label{cor2}
The cohomology class of the  
central extension $\mathbf{\widetilde{\Gamma^s_{g,r}}}$ is 
\[c_{{\mathbf{\widetilde{\Gamma^s_{g,r}}}}}= 12 \chi +\sum_{i=1}^s e_i \in 
H^2(\Gamma^s_{g,r};\C^*)\] 
if $g\geq 3$ and $s\geq 4$. The same formula holds 
also when $g=2$ but the class $\chi$ vanishes in 
$H^2(\Gamma^s_{g,r};\C^*)$.  
Here $\chi$ and $e_i$ are one fourth of the Meyer signature class and 
respectively the $i$-th Euler class with $\C^*$ coefficients.
\end{corollary}

\begin{remark}
The central extension arising from
$SU(2)$-TQFT with $p_1$-structures
was computed in \cite{Ge,MR} for $\Gamma_{g}$ and it equals $12\chi$. 
It can be shown that their computations extend to the case of 
punctured surfaces and the associated class  for $\Gamma_{g,r}^s$ is 
$12 \chi +\sum_{i=1}^s e_i$.  Our result shows that this extension 
coincides with the central extension arising from quantum Teichm\"uller theory.
\end{remark}

\vspace{0.1cm}\noindent
The organization of the paper is as follows. In Section~\ref{sec:1}, we review the quantization of the
Teichm\"uller space of a punctured surface and define
the associated quantum representations of the decorated Ptolemy groupoid which correspond to linear
representations of a central extension of the decorated Ptolemy groupoid. Then, in Section~\ref{sec:2}, we prove Theorem~\ref{dilog} by finding the pull-back of this central
extension to the mapping class group of the surface. The key idea is to use a
Grothendieck type principle. Namely, one can identify
a central extension of the mapping class group of some
surface, if one understands its restrictions
to the mapping class groups of sub-surfaces of
bounded topological types.
The core of the proof consists in computing explicitly
the lifts to the central extension of the decorated Ptolemy groupoid of the relations known to hold in the mapping class groups.
When properly interpreted, these lifts yield the class of the
mapping class group extension.
\subsection*{Acknowledgements} The authors are indebted to
Stephane Baseilhac and Vlad Sergiescu
for useful discussions and to the referee for his/her
suggestions which helped improving the presentation.

\section{Quantum Teichm\"uller theory}\label{sec:1}
\subsection{The groupoid of decorated ideal triangulations}
Let $\Sigma=\Sigma_{g,r}^s$ be an oriented closed surface of genus $g$ with 
$r$ boundary components and $s\geq 1$ punctures. When $r>0$ we choose  
a set of points on each boundary, which will be called {\em boundary} 
punctures. When we need to single out the $s$ punctures lying in the interior 
we will call them {\em interior} punctures.
In this paper we will only consider 
the situation when each boundary component 
has exactly one boundary puncture, so that there is a total of $s+r$ 
punctures among which $r$ are boundary punctures. 
The triangulations of $\Sigma^s_{g,r}$  
whose vertices are the $s+r$ punctures  
will be called {\em ideal triangulations}. 
Then $\Sigma$ is {\em large} if and only if 
$\EULER s>0$, where $\EULER=4g-4+2s+3r$ is the number of 
triangles in an ideal triangulation.

\begin{definition}
A \emph{decorated ideal triangulation} of $\Sigma$ is an ideal triangulation
$\tau$ up to isotopy fixing the boundary, where all triangles are provided with a marked corner, and a 
bijective ordering map
    \[
\bar\tau\colon
\{1,\ldots,\EULER\}\ni j\mapsto\bar\tau_j\in\TRIANGLES(\tau)
    \]
is fixed. Here $\TRIANGLES(\tau)$ is the set of all triangles of $\tau$.
\end{definition}

\vspace{0.2cm}
\noindent
Graphically,  the marked corner of a triangle is indicated by an asterisk and the corresponding number is put inside the triangle. The set of all decorated ideal triangulations of  $\SURFACE$ is denoted
$\SDIT$.

\vspace{0.2cm}
\noindent
Recall that if a group $G$ freely acts on a set $X$, then there is an associated groupoid defined as follows. The objects are the $G$-orbits in $X$, while morphisms are $G$-orbits in $X\times X$ with respect to the diagonal action. Denote by $[x]$ the object represented by an element $x\in X$ and $[x,y]$ the morphism represented by a pair of elements $(x,y)\in X\times X$. Two morphisms $[x,y]$ and $[u,v]$, are composable if and only if $[y]=[u]$ and their composition is  $[x,y][u,v]=[x,gv]$, where $g\in G$ is the unique element sending $u$ to $y$. The inverse and the identity morphisms are given respectively by $[x,y]^{-1}=[y,x]$ and $\operatorname{id}_{[x]}=[x,x]$. In what follows,  products of the form $[x_1,x_2][x_2,x_3]\cdots[x_{n-1},x_n]$ will be shortened as $[x_1,x_2,x_3,\ldots,x_{n-1},x_n]$.

\vspace{0.2cm}
\noindent
The mapping class group $\Gamma^s_{g,r}$ of $\Sigma$ acts
freely on $\SDIT$. In this case, we let $\gpd$ denote the corresponding
groupoid, called the \emph{groupoid of decorated ideal triangulations},
or {\em decorated Ptolemy groupoid}. This  groupoid first considered 
in \cite{Ka} is an enhanced version of the usual Ptolemy groupoid 
introduced and studied by Penner in \cite{pe0} (see also \cite{pe1}), 
which arises in the Fock-Goncharov quantization  (\cite{FG1,FG}) 
of the Teichm\"uller space.
There is a presentation  for $\gpd$ with three types of
generators and four types of relations.

\vspace{0.2cm}
\noindent
The generators are of the form
$[\tau,\tau^\sigma]$, $[\tau,\rho_i\tau]$, and
$[\tau,\omega_{i,j}\tau]$, where $\tau^\sigma$
is obtained from $\tau$ by replacing the
ordering map $\bar\tau$ by the map
$\bar{\tau}\circ\sigma$, where
$\sigma\in S_{\EULER}$ is a permutation of
the set $\{1,\ldots,\EULER\}$, $\rho_i\tau$ is obtained from
$\tau$ by changing the marked corner of triangle
$\bar\tau_i$ as in Figure~\ref{ft}, and $\omega_{i,j}\tau$
is obtained from $\tau$ by applying the flip transformation
in the quadrilateral composed of triangles
$\bar\tau_i$ and $\bar\tau_j$ as in Figure~\ref{fe}.
\begin{figure}[htb]
  \centering
\begin{picture}(200,20)
\put(0,0){\begin{picture}(40,20)
\put(0,0){\line(1,0){40}}
\put(0,0){\line(1,1){20}}
\put(20,20){\line(1,-1){20}}
\put(0,0){\circle*{3}}
\put(20,20){\circle*{3}}
\put(40,0){\circle*{3}}
\footnotesize
\put(33,0){$*$}
\put(18,5){$i$}
\end{picture}}
\put(160,0){\begin{picture}(40,20)
\put(0,0){\line(1,0){40}}
\put(0,0){\line(1,1){20}}
\put(20,20){\line(1,-1){20}}
\put(0,0){\circle*{3}}
\put(20,20){\circle*{3}}
\put(40,0){\circle*{3}}
\footnotesize
\put(17.5,14){$*$}
\put(18,5){$i$}
\end{picture}}
\put(95,8){$\stackrel{\rho_i}{\longrightarrow}$}
\end{picture}
\caption{The transformation $\rho_i$.}\label{ft}
\end{figure}
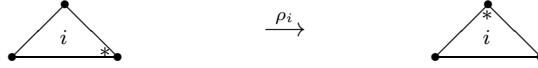
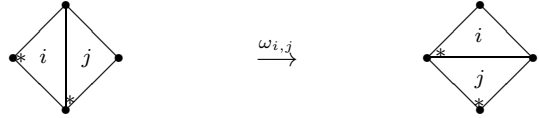
\begin{figure}[htb]
  \centering
\begin{picture}(200,40)
\put(0,0){
\begin{picture}(40,40)
\put(20,0){\line(-1,1){20}}
\put(40,20){\line(-1,-1){20}}
\put(0,20){\line(1,1){20}}
\put(40,20){\line(-1,1){20}}
\put(20,0){\line(0,1){40}}
\put(20,0){\circle*{3}}
\put(0,20){\circle*{3}}
\put(20,40){\circle*{3}}
\put(40,20){\circle*{3}}
\footnotesize
\put(10,18){$i$}\put(26,18){$j$}
\put(1,18){$*$}
\put(19.5,2){$*$}
\end{picture}}
\put(160,0){\begin{picture}(40,40)
\put(20,0){\line(-1,1){20}}
\put(40,20){\line(-1,-1){20}}
\put(0,20){\line(1,1){20}}
\put(40,20){\line(-1,1){20}}
\put(0,20){\line(1,0){40}}
\put(20,0){\circle*{3}}
\put(0,20){\circle*{3}}
\put(20,40){\circle*{3}}
\put(40,20){\circle*{3}}
\footnotesize
\put(18,26){$i$}\put(18,10){$j$}
\put(3,20){$*$}
\put(17.5,1){$*$}
\end{picture}}
\put(95,17){$\stackrel{\omega_{i,j}}{\longrightarrow}$}
\end{picture}
\caption{The transformation $\omega_{i,j}$.}\label{fe}
\end{figure}

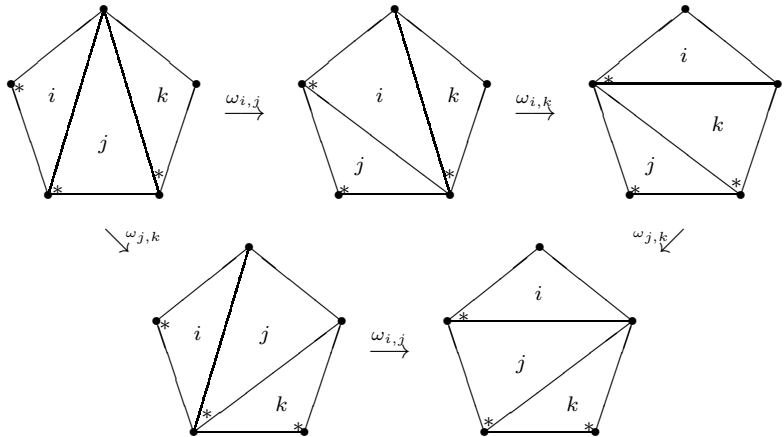
\begin{figure}[htb]
  \centering
\begin{picture}(290,160)
\put(0,90){\begin{picture}(70,70)
\put(14,0){\line(-1,3){14}}
\put(56,0){\line(1,3){14}}
\put(0,42){\line(5,4){35}}
\put(70,42){\line(-5,4){35}}
\put(14,0){\line(1,0){42}}
\qbezier(14,0)(20,20)(35,70)
\qbezier(56,0)(50,20)(35,70)
\put(14,0){\circle*{3}}
\put(56,0){\circle*{3}}
\put(0,42){\circle*{3}}
\put(70,42){\circle*{3}}
\put(35,70){\circle*{3}}
\footnotesize
\put(1,38.5){$*$}
\put(15.5,0){$*$}
\put(53.75,6){$*$}
\put(14,35){$i$}
\put(33,17){$j$}
\put(55,35){$k$}
\end{picture}}
\put(110,90){\begin{picture}(70,70)
\put(14,0){\line(-1,3){14}}
\put(56,0){\line(1,3){14}}
\put(0,42){\line(5,4){35}}
\put(70,42){\line(-5,4){35}}
\put(14,0){\line(1,0){42}}
\put(56,0){\line(-4,3){56}}
\qbezier(56,0)(50,20)(35,70)
\put(14,0){\circle*{3}}
\put(56,0){\circle*{3}}
\put(0,42){\circle*{3}}
\put(70,42){\circle*{3}}
\put(35,70){\circle*{3}}
\footnotesize
\put(2,39.5){$*$}
\put(14,0){$*$}
\put(53.75,6){$*$}
\put(28,35){$i$}
\put(20,9){$j$}
\put(55,35){$k$}
\end{picture}}
\put(220,90){\begin{picture}(70,70)
\put(14,0){\line(-1,3){14}}
\put(56,0){\line(1,3){14}}
\put(0,42){\line(5,4){35}}
\put(70,42){\line(-5,4){35}}
\put(14,0){\line(1,0){42}}
\put(56,0){\line(-4,3){56}}
\put(0,42){\line(1,0){70}}
\put(14,0){\circle*{3}}
\put(56,0){\circle*{3}}
\put(0,42){\circle*{3}}
\put(70,42){\circle*{3}}
\put(35,70){\circle*{3}}
\footnotesize
\put(4,41.5){$*$}
\put(14,0){$*$}
\put(52.25,2.5){$*$}
\put(33,50){$i$}
\put(20,9){$j$}
\put(45,24 ){$k$}
\end{picture}}
\put(55,0){\begin{picture}(70,70)
\put(14,0){\line(-1,3){14}}
\put(56,0){\line(1,3){14}}
\put(0,42){\line(5,4){35}}
\put(70,42){\line(-5,4){35}}
\put(14,0){\line(1,0){42}}
\qbezier(14,0)(20,20)(35,70)
\put(14,0){\line(4,3){56}}
\put(14,0){\circle*{3}}
\put(56,0){\circle*{3}}
\put(0,42){\circle*{3}}
\put(70,42){\circle*{3}}
\put(35,70){\circle*{3}}
\footnotesize
\put(1,38.5){$*$}
\put(17,5){$*$}
\put(51.5,0){$*$}
\put(14,35){$i$}
\put(39,35){$j$}
\put(45,8){$k$}
\end{picture}}
\put(165,0){\begin{picture}(70,70)
\put(14,0){\line(-1,3){14}}
\put(56,0){\line(1,3){14}}
\put(0,42){\line(5,4){35}}
\put(70,42){\line(-5,4){35}}
\put(14,0){\line(1,0){42}}
\put(0,42){\line(1,0){70}}
\put(14,0){\line(4,3){56}}
\put(14,0){\circle*{3}}
\put(56,0){\circle*{3}}
\put(0,42){\circle*{3}}
\put(70,42){\circle*{3}}
\put(35,70){\circle*{3}}
\footnotesize
\put(4,41.5){$*$}
\put(13.5,2){$*$}
\put(51.5,0){$*$}
\put(33,50){$i$}
\put(26,24){$j$}
\put(45,8){$k$}
\end{picture}}
\put(35,70){$\searrow$}\put(43,74){\tiny$\omega_{j,k}$}
\put(245,70){$\swarrow$}\put(235,74){\tiny$\omega_{j,k}$}
\put(135,28){$\stackrel{\omega_{i,j}}{\longrightarrow}$}
\put(80,118){$\stackrel{\omega_{i,j}}{\longrightarrow}$}
\put(190,118){$\stackrel{\omega_{i,k}}{\longrightarrow}$}
\end{picture}
  \caption{The Pentagon relation~\eqref{eq:22}.}
  \label{fig:pen-om}
\end{figure}

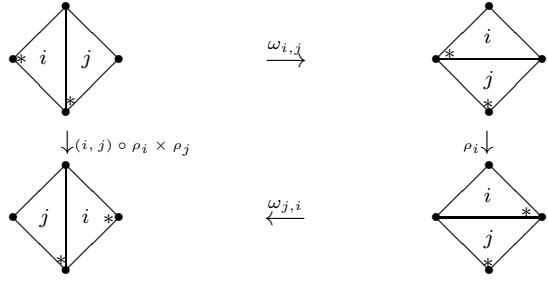
\begin{figure}[htb]
  \centering
  \begin{picture}(200,100)
\put(0,60){\begin{picture}(200,40)
\put(0,0){\begin{picture}(40,40)
\put(20,0){\line(-1,1){20}}
\put(40,20){\line(-1,-1){20}}
\put(0,20){\line(1,1){20}}
\put(40,20){\line(-1,1){20}}
\put(20,0){\line(0,1){40}}
\put(20,0){\circle*{3}}
\put(0,20){\circle*{3}}
\put(20,40){\circle*{3}}
\put(40,20){\circle*{3}}
\footnotesize
\put(10,18){$i$}\put(26,18){$j$}
\put(1,18){$*$}
\put(19.5,2){$*$}
\end{picture}}
\put(160,0){\begin{picture}(40,40)
\put(20,0){\line(-1,1){20}}
\put(40,20){\line(-1,-1){20}}
\put(0,20){\line(1,1){20}}
\put(40,20){\line(-1,1){20}}
\put(0,20){\line(1,0){40}}
\put(20,0){\circle*{3}}
\put(0,20){\circle*{3}}
\put(20,40){\circle*{3}}
\put(40,20){\circle*{3}}
\footnotesize
\put(18,26){$i$}\put(18,10){$j$}
\put(3,20){$*$}
\put(17.5,1){$*$}
\end{picture}}
\put(95,17){$\stackrel{\omega_{i,j}}{\longrightarrow}$}
\end{picture}}
\put(18,46){$\downarrow$\tiny$(i,j)\circ\rho_i\times\rho_j$}
\put(170,46){{\tiny$\rho_i$}$\downarrow$}
\put(0,0){\begin{picture}(200,40)
\put(0,0){\begin{picture}(40,40)
\put(20,0){\line(-1,1){20}}
\put(40,20){\line(-1,-1){20}}
\put(0,20){\line(1,1){20}}
\put(40,20){\line(-1,1){20}}
\put(20,0){\line(0,1){40}}
\put(20,0){\circle*{3}}
\put(0,20){\circle*{3}}
\put(20,40){\circle*{3}}
\put(40,20){\circle*{3}}
\footnotesize
\put(10,18){$j$}\put(26,18){$i$}
\put(16,2){$*$}
\put(34,17.5){$*$}
\end{picture}}
\put(160,0){\begin{picture}(40,40)
\put(20,0){\line(-1,1){20}}
\put(40,20){\line(-1,-1){20}}
\put(0,20){\line(1,1){20}}
\put(40,20){\line(-1,1){20}}
\put(0,20){\line(1,0){40}}
\put(20,0){\circle*{3}}
\put(0,20){\circle*{3}}
\put(20,40){\circle*{3}}
\put(40,20){\circle*{3}}
\footnotesize
\put(18,26){$i$}\put(18,10){$j$}
\put(32,20){$*$}
\put(17.5,1){$*$}
\end{picture}}
\put(95,17){$\stackrel{\omega_{j,i}}{\longleftarrow}$}
\end{picture}}
\end{picture}
  \caption{The Inversion relation~\eqref{eq:20}.}
  \label{fig:inv-rel}
\end{figure}

\vspace{0.2cm}
\noindent
There are two sets of relations satisfied by these generators. The first set is as follows:
\begin{gather}
\label{eq:23}
[\tau,\tau^\alpha,(\tau^\alpha)^\beta]=[\tau,\tau^{\alpha\beta}],\quad \alpha,\beta\in S_{\EULER},\\
  \label{eq:19}
  [\tau,\rho_i\tau,\rho_i\rho_i\tau,\rho_i\rho_i\rho_i\tau]=\operatorname{id}_{[\tau]},\\\label{eq:22}
  [\tau,\omega_{ij}\tau,\omega_{ik}\omega_{ij}\tau,\omega_{jk}\omega_{ik}\omega_{ij}\tau]
  =[\tau,\omega_{jk}\tau,\omega_{ij}\omega_{jk}\tau],
\\
\label{eq:20}
[\tau,\omega_{ij}\tau,\rho_i\omega_{ij}\tau,\omega_{ji}\rho_i\omega_{ij}\tau]
=[\tau,\tau^{(ij)},\rho_j\tau^{(ij)},\rho_i\rho_j\tau^{(ij)}],
\end{gather}
where the first two relations are evident, while the other two are shown graphically in Figures~\ref{fig:pen-om},~\ref{fig:inv-rel}.

\noindent The following commutation relations fulfill the second set of relations:
\begin{gather}
[\tau,\rho_i\tau,(\rho_i\tau)^\sigma]=[\tau,\tau^\sigma,\rho_{\sigma^{-1}(i)}\tau^\sigma],\\
[\tau,\omega_{ij}\tau,(\omega_{ij}\tau)^\sigma]=[\tau,\tau^\sigma,
\omega_{\sigma^{-1}(i)\sigma^{-1}(i)}\tau^\sigma],\\
[\tau,\rho_j\tau,\rho_i\rho_j\tau]=[\tau,\rho_i\tau,\rho_j\rho_i\tau],\\
[\tau,\rho_i\tau,\omega_{jk}\rho_i\tau]=[\tau,\omega_{jk}\tau,\rho_i\omega_{jk}\tau],\ i\not\in\{j,k\},\\
[\tau,\omega_{ij}\tau,\omega_{kl}\omega_{ij}\tau]=[\tau,\omega_{kl}\tau,\omega_{ij}\omega_{kl}\tau],\
\{i,j\}\cap\{k,l\}=\emptyset.
\end{gather}

\vspace{0.2cm}\noindent 
Consider now an embedding of $\Sigma^s_{g,r}$ into $\Sigma^{t}_{h,v}$ 
sending all punctures (both interior and boundary) 
to punctures. Of course boundary punctures are sent into interior punctures
unless the respective boundary circle is also a boundary of the larger 
surface. 

\begin{lemma}\label{embedgroupoid}
Assume that each component of $\Sigma^{t}_{h,v}\setminus {\rm int}(\Sigma^s_{g,r})$ 
is large. Then there is a natural embedding of 
$\mathcal G_{\Sigma^s_{g,r}}$ into  $\mathcal G_{\Sigma^t_{h,v}}$.  
\end{lemma}
\begin{proof}
Let $\tau_{ext}$ be a fixed decorated triangulation of 
$\Sigma^{t}_{h,v}\setminus {\rm int}(\Sigma^s_{g,r})$. If $\tau$ is a 
decorated triangulation of $\Sigma^s_{g,r}$ we denote by 
$\tau\cup \tau_{ext}$ the result of gluing the two triangulations 
along their corresponding boundary circles with the induced decoration. 
The isotopy class of the resulting triangulation is unique up to 
the action of Dehn twists along boundary components of $\Sigma^s_{g,r}$. 
This induces an injective map 
between the set of objects of the two groupoids. 
Then, the map which associates 
to the class $[\tau_1,\tau_2]$ of decorated triangulations of $\Sigma^s_{g,r}$ 
the class $[\tau_1\cup\tau_{ext},\tau_2\cup\tau_{ext}]$ is well-defined.  
Since the restriction of a homeomorphism of 
$\Sigma^{t}_{h,v}$ preserving the isotopy class 
of the decorated triangulation $\tau_{ext}$  to 
$\Sigma^{t}_{h,v}\setminus {\rm int}(\Sigma^s_{g,r})$ is isotopic to identity 
by Alexander's trick, the map defined above is injective.   
\end{proof}

\begin{remark}
When $r>0$ the construction of the decorated Ptolemy groupoid 
$\mathcal G_{\Sigma^s_{g,r}}$ 
depends on the choice of the set of boundary punctures, 
which might have more than $r$ elements, in general.  
\end{remark}

\subsection{Hilbert spaces of square integrable functions associated to triangulations}
\label{sec:q-functor}

In what follows, we work with Hilbert spaces
\[
\vsp\equiv L^2(\REALS),\quad \vsp^{\otimes n}\equiv L^2(\REALS^n).
\]
Any two self-adjoint operators
$\MOM$ and $\POSI$, acting in $\vsp$ and satisfying the Heisenberg commutation relation
\begin{equation}\label{E:heisen}
\MOM\POSI-\POSI\MOM=(2\pi \IMUN)^{-1}\operatorname{id}_{\vsp},
\end{equation}
can be realized as differentiation and multiplication operators. Such "coordinate" realization in Dirac's bra-ket notation has the form
\begin{equation}\label{eq:mdo}
\langle x|\MOM =\frac1{2\pi\IMUN}\frac{\partial}{\partial x}\langle x|,
\quad \langle x|\POSI =x\langle x|.
\end{equation}
Formally, the set of "vectors" $\{|x\rangle\}_{x\in\mathbb{R}}$ forms a generalized basis of $\vsp$ with the following orthogonality and completeness properties:
\[
\langle x|y\rangle=\delta(x-y),\quad \int_{\mathbb{R}}|x\rangle dx\langle x|=\operatorname{id}_{\vsp}.
\]
For any $1\le i\le m$ we shall use the following notation
        \[
        \iota_i\colon \End\vsp\ni \sfa\mapsto
        \sfa_i=\underbrace{1\otimes\cdots\otimes1}_{i-1\ \mathrm{times}}\otimes
        \sfa\otimes1\otimes\cdots\otimes1
        \in\End\vsp^{\otimes m}.
         \]
Besides that, if $\sfu\in \End\vsp^{\otimes k}$ for some $1\le
         k\le m$ and  $\{i_1,i_2,\ldots,i_k\}\subset\{1,2,\ldots,m\}$,
         then we shall write
    \[
\sfu_{i_1i_2\ldots i_2}\equiv\iota_{i_1}\otimes\iota_{i_2}\otimes\cdots
\otimes\iota_{i_k}(\sfu).
    \]
The symmetric group $\PGROUP_m$ naturally acts in $\vsp^{\otimes m}$:
    \begin{equation}\label{eq:perm}
\PERMUTE_\sigma (x_1\otimes\cdots\otimes
x_i\otimes\cdots\otimes x_m) =x_{\sigma^{-1}(1)}\otimes\cdots\otimes
x_{\sigma^{-1}(i)}\otimes\ldots\otimes x_{\sigma^{-1}(m)},\quad\sigma\in \PGROUP_m.
    \end{equation}

\subsection{Semi-symmetric $T$-matrices}
\label{sec:inii-aeaa-neno}

We define now the algebraic structure needed for constructing representations
of the decorated Ptolemy groupoid $\gpd$.

\begin{definition}
A \emph{semi-symmetric $T$-matrix} consists of two operators
$\ROTATE\in \operatorname{End}(\vsp)$ and $\PTOLEMY \in \operatorname{End}(\vsp^{\otimes2})$ satisfying
the equations:
\begin{gather}
  \label{eq:1}
\ROTATE^3=1, \\
 \label{eq:2}
\PTOLEMY_{12}\PTOLEMY_{13}\PTOLEMY_{23}=\PTOLEMY_{23}\PTOLEMY_{12}, \\
  \label{eq:4}
\PTOLEMY_{12} \ROTATE_1\PTOLEMY_{21}=\zeta\ROTATE_1\ROTATE_2\PERMUTE_{(12)},
\end{gather}
where $\zeta\in \C^*$ and
the permutation operator $\PERMUTE_{(12)}$ is defined by
equation~\eqref{eq:perm},  for $\sigma$ denoting the transposition $(12)$.
\end{definition}

\vspace{0.2cm}\noindent
Examples of semi-symmetric $T$-matrices could be obtained as follows.
Fix some self-conjugate operators $\MOM,\POSI$ satisfying the
Heisenberg commutation relation\eqref{E:heisen}.
Choose a parameter $\la$ satisfying the condition:
\[
(1-|\la|){\rm Im}\la=0,
\]
and define then two unitary operators by the following formulas:
\begin{gather}
  \label{eq:5}
  \ROTATE\equiv e^{-\IMUN\pi/3}e^{\IMUN
    3\pi\POSI^2}e^{\IMUN\pi(\MOM+\POSI)^2}\in \End(\vsp),\\
\label{eq:t-in-t-of-psi}
\PTOLEMY\equiv e^{\IMUN 2\pi\MOM_1\POSI_2}
\QDILOG(\POSI_1+\MOM_2-\POSI_2)\in \End(\vsp^{\otimes2}).
    \end{gather}
They satisfy the defining relations for a semi-symmetric $T$-matrix,
where
\begin{equation}
  \label{eq:prfac}
  \zeta =  e^{\IMUN\pi\cla^2/3},\quad \cla=\frac\imun2(\la+\la^{-1}),
\end{equation}
and $\QDILOG$ is Faddeev's  quantum dilogarithm defined 
on $\{z\in \C; |\rm Im(z)| < |Im(\cla)|\}$ by means of 
\begin{equation}
\QDILOG(z)=\exp\left(-\frac{1}{4}\int_{-\infty}^{\infty}
\frac{\exp(-2 i zx) d\,x}{\sinh(xb)\sinh(x/b)x}\right)
\end{equation} 
Faddeev's quantum dilogarithm is closely related to the double gamma and double 
sine  functions (\cite{Bar,Shi}) and was used by  Baxter (\cite{Bax}) 
and Faddeev (see \cite{Fad,FaK}). 
Its main feature is the following functional equation  (see \cite{Fad,FaK}) 
it satisfies: 
\[ \QDILOG(q)\QDILOG(p)=\QDILOG(p)\QDILOG(p+q)\QDILOG(q)\]
whenever $pq-qp=\frac{1}{2\pi i}\bf 1$. 
 
\vspace{0.1cm}\noindent
Remark that  the operator $\ROTATE$ is
characterized (up to a normalization factor) by the equations:
\[
\ROTATE\POSI\ROTATE^{-1}=\MOM-\POSI,\quad\ROTATE\MOM\ROTATE^{-1}=-\POSI.
\]
Note that equations~(\ref{eq:1})---(\ref{eq:4}) correspond to relations~
\eqref{eq:19}---\eqref{eq:20}.


\vspace{0.2cm}\noindent
Let us introduce now some notation which will be useful in the sequel.
For any operator $\sfa\in\End \vsp$ we set:
\begin{equation}
  \label{eq:14}
  \sfa_{\hat k}\equiv \ROTATE_k\sfa_k\ROTATE_k^{-1},\quad
 \sfa_{\check k}\equiv \ROTATE_k^{-1}\sfa_k\ROTATE_k.
\end{equation}
It is evident that
\[
\sfa_{\check{\hat k}}=\sfa_{\hat{\check k}}=\sfa_k,\quad \sfa_{\hat{\hat
    k}}=\sfa_{\check k},\quad \sfa_{\check{\check k}}=\sfa_{\hat k},
\]
where the last two equations follow from equation~\eqref{eq:1}. In particular,
we have
\begin{gather}
\MOM_{\hat k}=-\POSI_k,\quad\POSI_{\hat k}=\MOM_k-\POSI_k,\\
\MOM_{\check k}=\POSI_k-\MOM_k,\quad\POSI_{\check k}=-\MOM_k.
\end{gather}
Besides that, it will be also useful to use the notation
\begin{equation}
  \label{eq:16}
 \PERMUTE_{(kl\ldots m\hat k)}\equiv\ROTATE_k\PERMUTE_{(kl\ldots
  m)},\quad
\PERMUTE_{(kl\ldots m\check k)}\equiv\ROTATE_k^{-1}\PERMUTE_{(kl\ldots
  m)},
\end{equation}
where $(kl\ldots m)$ is the cyclic permutation
\[
(kl\ldots m)\colon k\mapsto l\mapsto\ldots\mapsto m\mapsto k.
\]
Equation~\eqref{eq:4} in this notation takes a rather compact form
\begin{gather}
\label{eq:16a}
  \PTOLEMY_{12}\PTOLEMY_{2\hat1}= \zeta\PERMUTE_{(12\hat1)}.
\end{gather}
\begin{remark}
Notice that the Pentagon relation (\ref{eq:2}) 
can be applied whenever any of the  
indices $k\in\{1,2\}$ arising among subscripts 
is replaced everywhere by either $\hat k$ or else 
$\check k$. 
\end{remark}
\begin{remark}\label{symmetry}
A $T$-matrix has the following symmetry property:
$\PTOLEMY_{12}=\PTOLEMY_{\hat2\check1}$. 
This can be obtained using twice relation (\ref{eq:16a}): 
\begin{equation}\label{symm}
T_{12}=T_{12}T_{2\hat1}T_{2\hat1}^{-1}=\zeta P_{(12\hat1)}T_{2\hat1}^{-1}=
T_{\hat1\hat2}^{-1}\zeta P_{(12\hat1)}=
T_{\hat1\hat2}^{-1}\zeta P_{(\hat1\hat2\check1)}=T_{\hat2\check1}
\end{equation}
\end{remark}

\subsection{The quantum Teichm\"uller space}
\label{sec:eaaioiaue-ooieoid}
The quantization of the Teichm\"uller space of a punctured surface $\Sigma$
with boundary induced by a semi-symmetric $T$-matrix is defined by means of a
\emph{quantum functor}:
\[
\FUNCTOR\colon\gpd\to\End(\vsp^{\otimes\NTR}),
\]
Its meaning is that we have an operator valued function:
\[
\FUNCTOR\colon\SDIT\times\SDIT\to \End(\vsp^{\otimes\NTR}),
\]
satisfying the  following equations:
 \begin{equation}
  \label{eq:10}
 \FUNCTOR(\tau,\tau)=\operatorname{id}_{\vsp^{\otimes\NTR}},\quad \FUNCTOR(\tau,\tau')\FUNCTOR(\tau',\tau'')
\FUNCTOR(\tau'',\tau)\in\COMPLEXS\setminus\{0\},\quad \forall\tau,\tau',\tau''\in\SDIT,
\end{equation}
 \begin{equation}
\FUNCTOR(f(\tau),f(\tau'))= \FUNCTOR(\tau,\tau'),\qquad \forall
f\in\MCG,
\end{equation}
\begin{equation}
  \label{eq:17}
     \FUNCTOR(\tau,\rho_i\tau)\equiv\ROTATE_i,
\end{equation}
    \begin{equation}\label{tij}
\FUNCTOR(\tau,\omega_{i,j}\tau)\equiv\PTOLEMY_{ij},
    \end{equation}
\begin{equation}
\FUNCTOR(\tau,\tau^\sigma)\equiv\PERMUTE_\sigma,
\quad\forall\sigma\in\PGROUP_{\NTR},
    \end{equation}
where operator $\PERMUTE_\sigma$ is defined by equation~(\ref{eq:perm}).
Consistency of these equations is ensured by the consistency of equations~(\ref{eq:1})---(\ref{eq:4}) with relations~\eqref{eq:19}---\eqref{eq:20}.

\vspace{0.2cm}
\noindent
A particular case of equation~(\ref{eq:10}) corresponds to
$\tau''=\tau$:
\begin{equation}
  \label{eq:18}
  \FUNCTOR(\tau,\tau')\FUNCTOR(\tau',\tau)\in\COMPLEXS\setminus\{0\}.
\end{equation}
As an example, we can calculate the operator
$\FUNCTOR(\tau,\omega_{i,j}^{-1}(\tau))$. Denoting $\tau'\equiv
\omega_{i,j}^{-1}(\tau)$ and using equation~(\ref{eq:18}), as well as definition~\eqref{tij}, we obtain
\begin{equation}
  \label{eq:11}
 \FUNCTOR(\tau,\omega_{i,j}^{-1}(\tau))=
\FUNCTOR(\omega_{i,j}(\tau'),\tau')\simeq
(\FUNCTOR(\tau',\omega_{i,j}(\tau')))^{-1}=\PTOLEMY_{ij}^{-1},
\end{equation}
where  $\simeq$ means equality up to a numerical multiplicative factor.

\vspace{0.2cm}
\noindent 
The operations $\hat{}$ and $\check{}$ at the indices level 
have the following geometric interpretation.  
If the distinguished corners of the decorated ideal triangulation 
are precisely those from Figure \ref{fe} 
then the quantum functor assigns to the flip on that edge 
the endomorphism $T_{ij}^{-1}$. Now, changing the distinguished 
corner in the triangle labeled $i$ amounts of changing $i$ into 
$\hat i$ or $\check i$ (and similarly for $j$) in the expression 
of the quantum functor endomorphism. This rules will be 
intensively used when we compute the expressions of 
Dehn twists in terms of the generators of the 
decorated Ptolemy groupoid in the next section. 

\vspace{0.2cm}
\noindent
The quantum functor induces a unitary projective  representation
of the mapping class group $\Gamma^s_g$ of $\Sigma$ as follows:
    \[
    \Gamma^s_g\ni f\mapsto\FUNCTOR(\tau,f(\tau))\in
\End (\vsp^{\otimes\NTR}).
    \]
 Indeed, we have the following relation (up to a non-zero scalar):
\[
\FUNCTOR(\tau,f(\tau))\FUNCTOR(\tau,h(\tau))=
\FUNCTOR(\tau,f(\tau))\FUNCTOR(f(\tau),f(h(\tau)))
\simeq\FUNCTOR(\tau,fh(\tau)).
\]
The main question addressed in this present paper is to identify the central extension
of the mapping class group corresponding to this projective representation.
Observe that the projective factor lies in the sub-group of $\mathbb{C}^*$ generated by $\zeta$.

\vspace{0.2cm}
\noindent 
In \cite{Ka,Ka2} one considered only punctured surfaces 
without boundary. However, the construction extends 
without essential modifications to the case when $\Sigma$ 
is  a surface with boundary $\Sigma^s_{g,r}$ when 
$s\geq 1$ and each boundary component contains  
one boundary puncture. In this case we could define directly 
the central extension $\widetilde{\Gamma^s_{g,r}}$ by using the 
decorated Ptolemy groupoid of the punctured surface with boundary, 
without reference to a larger surface without boundary.

\section{Presentation of  $\widetilde{\Gamma^s_{g,r}}$}\label{sec:2}
\subsection{Generating set for the relations}
We start with a number of notations and definitions.
Our setup consists of an embedding $\Sigma^s_{g,r}\subset  
\Sigma^{t}_{h,0}$ sending punctures into punctures.
We assume that each component of 
$\Sigma^{t}_{h,0}\setminus {\rm int}(\Sigma^s_{g,r})$ is large, 
namely it admits ideal triangulations whose vertices are  
those punctures of $\Sigma^{t}_{h,0}$ which are not interior 
punctures of $\Sigma^s_{g,r}$ (hence boundary punctures of 
$\Sigma^s_{g,r}$ being allowed). 
In particular,  if we discard the boundary punctures of 
$\Sigma^s_{g,r}$ the complement 
$\Sigma^{t}_{h,0}\setminus {\rm int}(\Sigma^s_{g,r})$ contains no disk, 
punctured disk or cylinder components. 
According to \cite{PR1} the surface embedding 
induces an embedding between the corresponding mapping class groups 
 $\Gamma^s_{g,r}\hookrightarrow  \Gamma^{t}_{h,0}$. 
The pull-back of the central extension 
$\widetilde{\Gamma^t_{h}}$ to $\Gamma^s_{g,r}$ is a central extension 
$\widetilde{\Gamma^s_{g,r}}$. Our main concern is to study this 
central extension. The central extension obtained by 
the present construction is isomorphic to the central extension 
obtained by the direct quantization of the Teichm\"uller space 
associated to $\Sigma^s_{g,r}$ following the procedure of section 1.4.  
This follows from the fact that the map between the  
mapping class groups $\Gamma^s_{g,r}\hookrightarrow  \Gamma^{t}_{h,0}$ 
is covered by an injective map between the decorated Ptolemy groupoids 
according to Lemma \ref{embedgroupoid}.

\vspace{0.2cm}\noindent 
Since the restriction of the 
Euler class corresponding to the $(s+1)$-th puncture to $\Gamma^s_{g,r}$ 
vanishes, it is enough to consider $t=s$ below. 
Our strategy is to compute explicit lifts to  
$\widetilde{\Gamma^s_{g,r}}$ of a set of relations arising in a 
group presentation of  $\Gamma^s_{g,r}$ by expressing (lifts of) the 
generators as elements of the decorated Ptolemy groupoid of the 
larger punctured surface $\Sigma^{s}_{h,0}$.  
The independence on the particular embedding of the subsurface 
$\Sigma^s_{g,r}$, under the assumptions 
of the main theorem is a consequence of the so-called Grothendieck principle.  
In the form proved by Gervais in \cite{Ge} it 
states that all relations in $\Gamma^s_{g,r}$ are determined 
by an explicit set of relations among mapping classes supported on 
small subsurfaces, namely $\Sigma_{0,4}$, $\Sigma_{1,2}$ and 
$\Sigma_{0,3}$, where $\Sigma_{g,r}=\Sigma^0_{g,r}$. 
We express then these relations in terms of elements of the decorated 
Ptolemy groupoids of the surfaces $\Sigma^4_{0,4}$, $\Sigma^2_{1,2}$ and 
$\Sigma^4_{0,3}$, respectively. According to Lemma \ref{embedgroupoid} 
these relations also hold in $\mathcal G_{\Sigma^s_{h,0}}$, 
provided that $s\geq 4$.

\vspace{0.2cm}\noindent   
If $a$ is a simple closed curve on $\Sigma^s_{g,r}$ 
we denote by $D_a\in \Gamma^s_{g,r}$ the {\em right} Dehn twist along $a$. 

\begin{definition}
A {\em chain relation} $C$ on the surface $\Sigma^s_{g,r}$
is given by an embedding $\Sigma_{1,2}\subset \Sigma^s_{g,r}$
and the  standard chain relation on this 2-holed torus, namely
\[ (D_aD_bD_c)^4=D_eD_f\]
where $a,b,c,d,e,f$ are the following curves of the embedded 2-holed torus:

\vspace{0.2cm}
\begin{center}
\includegraphics[scale=0.4]{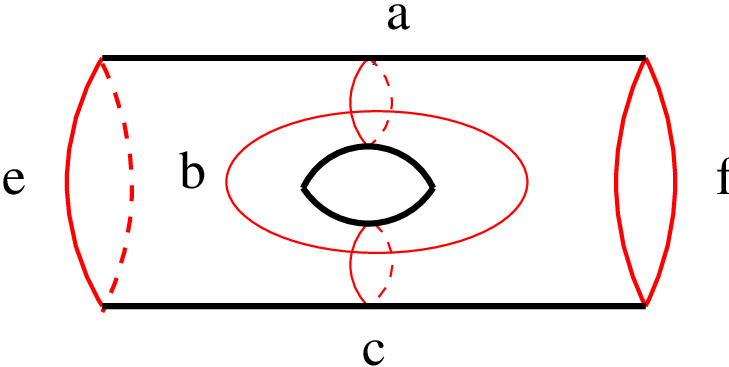}
\end{center}
\end{definition}

\begin{definition}
A {\em lantern relation} $L$ on the surface $\Sigma^s_{g,r}$
is given by an embedding $\Sigma_{0,4}\subset \Sigma^s_{g,r}$
and the  standard lantern relation on this 4-holed sphere, namely
\begin{equation} D_{a_{12}}D_{a_{13}}D_{a_{23}}D_{a_0}^{-1}D_{a_1}^{-1}D_{a_2}^{-1}D_{a_3}^{-1}=1
\end{equation}
where $a_0,a_1,a_2,a_3,a_{12},a_{13},a_{23}$
are the following curves of the embedded 4-holed sphere:

\vspace{0.2cm}
\begin{center}
\includegraphics[scale=0.4]{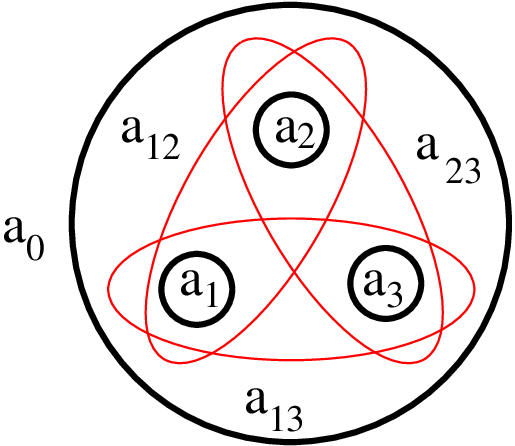}
\end{center}
\end{definition}

\begin{definition}
Consider an embedding $\Sigma^1_{0,3}\subset \Sigma^s_{g,r}$ such that the
boundary components $a_1,a_2,a_3$ of $\Sigma^1_{0,3}$ are non-separating curves.
Let then $a_{12},a_{13},a_{23}$ be embedded curves on $\Sigma^1_{0,3}$ so that
$a_{jk}$ bounds a pair of pants $\Sigma_{0,3}\subset \Sigma^1_{0,3}$ along with $a_j$ and $a_k$,
for all $1\leq j\neq k\leq 3$. Then the {\em puncture relation} $P$ (supported
at the puncture of $\Sigma^1_{0,3}$) on  the surface $\Sigma^s_{g,r}$ is:
\begin{equation}  
D_{a_{12}}D_{a_{13}}D_{a_{23}}D_{a_1}^{-1}D_{a_2}^{-1}D_{a_3}^{-1}=1\end{equation}
\end{definition}
\begin{remark}
The puncture relation is, in fact, a consequence of the
lantern relation  and the
fact that the Dehn twist along a small loop encircling a puncture is trivial.
\end{remark}

\vspace{0.2cm}\noindent
The first step in proving Theorem~\ref{dilog} is to find an explicit presentation
for the central extension $\widetilde{\Gamma^s_{g,r}}$. Specifically, by using Gervais' presentation~\cite{Ge}, we have the following description.
\begin{proposition}\label{pres}
Suppose that $g\geq 2$ and $s\geq 4$. Then 
the group $\widetilde{\Gamma^s_{g,r}}$ has the following presentation.
\begin{enumerate}
\item Generators:
\begin{enumerate}
\item With each non-separating  simple closed curve $a$
in $\Sigma^s_{g,r}$ is associated a generator  $\widetilde{D}_{a}$;
\item One (central) element $z$.
\end{enumerate}
\item Relations:
\begin{enumerate}
\item Centrality:
\begin{equation}
z \widetilde{D}_a=\widetilde{D}_az
\end{equation}
for any non-separating  simple closed curve $a$
on $\Sigma^s_{g,r}$;
\item Braid type $0$-relations:
\begin{equation}
\widetilde{D}_a\widetilde{D}_b=\widetilde{D}_b\widetilde{D}_a\end{equation}
for each pair of disjoint non-separating simple closed curves $a$ and $b$;
\item   Braid type $1$-relations:
\begin{equation} \widetilde{D}_a\widetilde{D}_b\widetilde{D}_a=\widetilde{D}_b\widetilde{D}_a
\widetilde{D}_b
\end{equation}
for each pair of non-separating simple closed curves $a$ and $b$ which intersect transversely at one point;
\item One lantern relation on a $4$-holed sphere
subsurface with non-separating boundary curves:
\begin{equation}
 \widetilde{D}_{a_0}\widetilde{D}_{a_1}\widetilde{D}_{a_2}\widetilde{D}_{a_3}=
\widetilde{D}_{a_{12}}\widetilde{D}_{a_{13}}\widetilde{D}_{a_{23}}
\end{equation}
\item One chain relation on a $2$-holed torus subsurface with non-separating
boundary curves:
\begin{equation} (\widetilde{D}_a\widetilde{D}_b\widetilde{D}_c)^4=z^{12}
\widetilde{D}_e\widetilde{D}_f \end{equation}
\item Puncture relations:
\begin{equation}
\widetilde{D}_{a_{12(i)}}\widetilde{D}_{a_{13}(i)}\widetilde{D}_{a_{23}(i)}
=z\widetilde{D}_{a_1(i)}\widetilde{D}_{a_2(i)}\widetilde{D}_{a_3(i)}
\end{equation}
for each puncture $p_i$ of $\Sigma^s_{g,r}$,  $i\in\{1,2,\ldots,s\}$.
\item Scalar equation: 
\begin{equation}
z^N=1
\end{equation} 
where $N$ is the order of $\zeta^{-6}$, in the case where $\zeta\in \C^*$
is a root of unity. 
\end{enumerate}
\end{enumerate}
\end{proposition}

\subsection{Proof of Proposition \ref{pres}}
\begin{lemma}
For any lifts  $\widetilde{D}_a$ of the Dehn twists $D_a$ we have
$\widetilde{D}_a\widetilde{D}_b=\widetilde{D}_b\widetilde{D}_a$, for any 
two disjoint simple closed curves $a$ and $b$, and
thus the braid-type 0-relations (b) are satisfied.
\end{lemma}
\begin{proof}
The commutativity relations are satisfied
for particular lifts coming from a semi-symmetric $T$-matrix.
If we change the lifts by multiplying each lift by some central element
the  commutativity is still valid. Thus, the commutativity holds for any lifts.
\end{proof}

\begin{lemma}\label{braid1}
There are lifts $\widetilde{D}_a$ of the Dehn twists $D_a$, for each
non-separating simple closed curve $a$  such that we have
$\widetilde{D}_a\widetilde{D}_b\widetilde{D}_a=\widetilde{D}_b\widetilde{D}_a\widetilde{D}_b$
for any simple closed curves $a,b$  with one intersection point, and thus
the braid type $1$-relations (c) are satisfied.
Moreover, the choice of lifts of all $\widetilde{D}_x$, with $x$ non-separating,
satisfying these requirements is uniquely defined by
fixing the lift $\widetilde{D}_{a}$ of one particular Dehn twist.
\end{lemma}
\begin{proof}
Consider an arbitrary lift of one braid type $1$-relation (to be called the fundamental one),
which has the form
$\widetilde{D}_a\widetilde{D}_b\widetilde{D}_a=z^k\widetilde{D}_b\widetilde{D}_a\widetilde{D}_b
$.   Change then  the lift $\widetilde{D}_b$ into $z^k\widetilde{D}_b$. With the new lift the
relation above becomes $\widetilde{D}_a\widetilde{D}_b\widetilde{D}_a=\widetilde{D}_b\widetilde{D}_a\widetilde{D}_b$.

\vspace{0.2cm}\noindent
Choose now an arbitrary braid type $1$-relation of $\Gamma^s_{g,r}$, say
$D_xD_yD_x=D_yD_xD_y$.   There exists a
1-holed torus $\Sigma_{1,1}\subset \Sigma^s_{g,r}$ containing $x,y$, namely a neighborhood of $x\cup y$.
Let $T$ be the similar torus containing $a,b$.  Since $a,b$ and $x,y$ are
non-separating there exists a
homeomorphism $\varphi:\Sigma^s_{g,r}\to \Sigma^s_{g,r}$ such that
$\varphi(a)=x$ and $\varphi(b)=y$. We have then
\[ D_x=\varphi D_a \varphi^{-1},\quad  D_y= \varphi D_b \varphi^{-1}.\]
Let us consider now an arbitrary lift $\widetilde{\varphi}$
of $\varphi$, which is well-defined only up to a central element, and set
\[ \widetilde{D}_x=\widetilde{\varphi} \widetilde{D}_a\widetilde{\varphi}^{-1},\quad
\widetilde{D}_y=\widetilde{\varphi} \widetilde{D}_b\widetilde{\varphi}^{-1}. \]
These lifts are well-defined since they do not depend on the choice of $\widetilde{\varphi}$ (the central elements  coming from $\widetilde{\varphi}$ and
$\widetilde{\varphi}^{-1}$ mutually cancel). Moreover, we have then
\[ \widetilde{D}_x\widetilde{D}_y\widetilde{D}_x=\widetilde{D}_y\widetilde{D}_x\widetilde{D}_y
\]
and so the braid type $1$-relations  (c) are all satisfied.

\vspace{0.2cm}\noindent
For the second part of the lemma observe that 
the choice of $\widetilde{D}_{a}$ fixes the choice of  $\widetilde{D}_{b}$.
If $x$ is a non-separating  simple closed curve on $\Sigma^s_{g,r}$, then
there exists another non-separating curve $y$ which intersects it in one point.
Thus, by the argument which was used above to prove the existence of the lifts 
the choice of $\widetilde{D}_{x}$ is unique.
\end{proof}

\begin{lemma}
One can choose the lifts of Dehn twists in
$\widetilde{\Gamma^s_{g,r}}$ so that all braid type relations are satisfied and
the lift of the lantern relation  (d) is
trivial, namely
 \[ 
 \widetilde{D}_{a_0}\widetilde{D}_{a_1}\widetilde{D}_{a_2}\widetilde{D}_{a_3}=
\widetilde{D}_{a_{12}}\widetilde{D}_{a_{13}}\widetilde{D}_{a_{23}}
\]
for the non-separating curves on an embedded $\Sigma_{0,4}\subset \Sigma^s_{g,r}$.
\end{lemma}
\begin{proof}
An arbitrary lift of that lantern relation is of the form
$\widetilde{D}_{a_0}\widetilde{D}_{a_1}\widetilde{D}_{a_2}\widetilde{D}_{a_3}=
z^k\widetilde{D}_{a_{12}}\widetilde{D}_{a_{13}}\widetilde{D}_{a_{23}}$.
In this case, we change the lift $\widetilde{D}_{a_0}$ into $z^{-k}\widetilde{D}_{a_0}$ and adjust the lifts of all other Dehn twists along non-separating  curves the way that all braid type $1$-relations are satisfied. Then, the required form of the lantern relation is satisfied.
\end{proof}

\vspace{0.2cm}\noindent
We say that the lifts of the Dehn twists are \emph{normalized} if all braid type
relations and one  lantern relation are lifted in a trivial way.
\begin{lemma}\label{normalization}
 Assume that $s\geq 4$. Then a normalized Dehn twist in quantum Teichm\"uller theory is conjugated to the inverse $T$-matrix times $\zeta^{-6}$ i.e.
\[
\widetilde{D}_\alpha=\fun(\tau,D_\alpha\tau)=\zeta^{-6}U_\alpha T_{kl}^{-1}U_\alpha^{-1}.
\]
\end{lemma}

\vspace{0.2cm}
\noindent 
As the computations  involved in the proof are rather laborious 
we postpone it after the proof of Lemma \ref{puncturel}. 

\vspace{0.2cm}\noindent 
We will suppose henceforth that the lifts of Dehn twists are normalized. 
\begin{lemma}\label{chainrel}
Let $a,b,c,e,f$ be the five curves appearing in the
chain relation \( (D_aD_bD_c)^4=D_eD_f\)on an embedded 2-holed torus sitting inside $\Sigma_{g,r}^s$.
If $s\geq 2$, then the lifts of Dehn twists in
$\widetilde{\Gamma^s_{g,r}}$ satisfy the relation
\[ (\widetilde{D}_a\widetilde{D}_b\widetilde{D}_c)^4=\zeta^{-72}\widetilde{D}_e\widetilde{D}_f  \]
\end{lemma}
\begin{proof}
\vspace{0.2cm}\noindent
If $s\geq 2$ and $g\geq 2$, then there is an embedding
$\Sigma^2_{2,1}\subset \Sigma^s_{g,r}$.

\vspace{0.2cm}
\noindent
We consider a surface $S$ homeomorphic to $\Sigma_{1,2}^2$, i.e. a torus with two holes and two punctures drawn in the left picture
of Figure~\ref{F:1} where the opposite sides of the
rectangle are identified.
Notice that the two punctures are located on the two boundary components.

\begin{figure}[h]
\begin{tikzpicture}[
]
\fill[blue!10] (-2,0) rectangle (2,2);
\draw[blue!50,fill=white] (0,2) .. controls (-1.5,1.5) and (-0.5,0.5) .. (0,2);
\draw[blue!50,fill=white] (2,2) .. controls (.5,1.5) and (1.5,0.5) .. (2,2);
\fill[white]
(-2,0) circle (1.5pt)
(0,0) circle (1.5pt)
(2,0) circle (1.5pt)
(-2,2) circle (1.5pt)
(0,2) circle (1.5pt)
(2,2) circle (1.5pt);
\end{tikzpicture}
\qquad
\begin{tikzpicture}[
]
\fill[blue!10] (-2,0) rectangle (2,2);
\draw[blue!50,fill=white] (0,2) .. controls (-1.5,1.5) and (-0.5,0.5) .. (0,2);

\draw[blue!50,fill=white] (2,2) .. controls (.5,1.5) and (1.5,0.5) .. (2,2);

\filldraw [white]
(-2,0) circle (1.5pt)
(0,0) circle (1.5pt)
(2,0) circle (1.5pt)
(-2,2) circle (1.5pt)
(0,2) circle (1.5pt)
(2,2) circle (1.5pt);
\draw[blue](-1.5,0)--(-1.5,2);
\draw[blue](.5,0)--(.5,2);
\draw[blue](-1,0)--(1,2);
\draw[blue] (1,0)--(2,1);
\draw [blue](-2,1)--(-1,2);
\draw[blue](0.2,0) arc (0:90:.2);
\draw[blue](0,0.2) arc (90:180:.5 and .2);
\draw[blue] (0.2,2) arc (0:-90:.2 and .5);
\draw[blue] (-.5,2) .. controls (-.5,1.8) and (-1.4,1.4).. (-1,1);
\draw[blue] (0,1.5) .. controls (-.2,1.5) and (-.6,.6).. (-1,1);

\draw[blue](-1.8,0) arc (0:90:.2);
\draw [blue](2,0.2) arc (90:180:.5 and .2);
\draw [blue](-1.8,2) arc (0:-90:.2 and .5);
\draw[blue](1.5,2) .. controls (1.5,1.8) and (.6,1.4).. (1,1);
\draw[blue](2,1.5) .. controls (1.8,1.5) and (1.4,.6).. (1,1);
\draw(-1.4,.5) node{\tiny $a$};
\draw(.6,.5) node{\tiny $c$};
\draw(-.3,.5) node{\tiny $b$};
\draw(-.9,1.1) node{\tiny $e$};
\draw(1.1,1.1) node{\tiny $f$};
\end{tikzpicture}
\qquad
\begin{tikzpicture}[
]

\draw[fill=blue!10] (-2,0) rectangle (2,2);
\draw[fill=white] (0,2) .. controls (-1.5,1.5) and (-0.5,0.5) .. (0,2);
\draw (0,2) .. controls (-.9,1.7) and (-1.7,.9) .. (-2,0);
\draw (0,2) .. controls (-.3,1.1) and (-1.1,.3) .. (-2,0);
\draw[fill=white] (2,2) .. controls (.5,1.5) and (1.5,0.5) .. (2,2);
\draw (2,2) .. controls (1.1,1.7) and (.3,.9) .. (0,0);
\draw (2,2) .. controls (1.7,1.1) and (.9,.3) .. (0,0);
\draw (0,2)--(0,0);
\draw (-1.9,.5) node {\tiny $1$};
\draw (-1.75,.25) node {\tiny $2$};
\draw (-1.4,.15) node {\tiny $3$};
\draw (.1,.5) node {\tiny $4$};
\draw (.25,.25) node {\tiny $5$};
\draw (.6,.15) node {\tiny $6$};
\filldraw [white]
(-2,0) circle (1.5pt)
(0,0) circle (1.5pt)
(2,0) circle (1.5pt)
(-2,2) circle (1.5pt)
(0,2) circle (1.5pt)
(2,2) circle (1.5pt);
\end{tikzpicture}
\caption{}\label{F:1}
\end{figure}
The central picture of Figure~\ref{F:1} specifies five simple closed curves $a,b,c,e,f$ in $S$, the Dehn twists along which enter the chain relation.

\vspace{0.2cm}
\noindent
 We also choose a particular decorated ideal triangulation $\tau$ of $S$ given by the right picture of Figure~\ref{F:1},
where the ideal arcs are drawn in black and the positions of the numbers in ideal triangles correspond to the marked corners. Notice that our choice is manifestly symmetric with respect to the exchange of the left and the right halves of the rectangle accompanied with relabeling $(1,2,3)\leftrightarrow (4,5,6)$. This symmetry will be useful for reducing the amount of calculations in deriving the quantum realizations of the Dehn twists.

\vspace{0.2cm}
\noindent
The basic procedure in deriving the quantum realization of the Dehn twist $D_\alpha$ along a given simple closed curve $\alpha$ is to use a specific decorated ideal triangulation where the contour $\alpha$ intersects only two ideal arcs, so that the annular neighborhood of $\alpha$ is given by only two ideal triangles. With respect to such (decorated) ideal triangulation the quantum operator realizing $D_\alpha$ is given by a single $T$-operator. Let us work out this procedure in the case of the curves $a,b,c,e,f$.

\vspace{0.2cm}
\noindent
For any simple closed curve $\alpha$, we denote $\bar\fun_\alpha=\widetilde{D}_\alpha^{-1}\simeq\fun(D_\alpha\tau,\tau)$.
To derive the operator representing the Dehn twist $D_a$, we apply the following change of triangulation:
\[
\xymatrix{
\begin{tikzpicture}[
baseline=(x.base)
]
\node (x) at (0,.8){};
\draw[fill=blue!10] (-2,0) rectangle (2,2);
\draw[fill=white] (0,2) .. controls (-1.5,1.5) and (-0.5,0.5) .. (0,2);
\draw (0,2) .. controls (-.9,1.7) and (-1.7,.9) .. (-2,0);
\draw (0,2) .. controls (-.3,1.1) and (-1.1,.3) .. (-2,0);
\draw[fill=white] (2,2) .. controls (.5,1.5) and (1.5,0.5) .. (2,2);
\draw (2,2) .. controls (1.1,1.7) and (.3,.9) .. (0,0);
\draw (2,2) .. controls (1.7,1.1) and (.9,.3) .. (0,0);
\draw (0,2)--(0,0);
\draw (-1.9,.5) node {\tiny $1$};
\draw (-1.75,.25) node {\tiny $2$};
\draw (-1.4,.15) node {\tiny $3$};
\draw (.1,.5) node {\tiny $4$};
\draw (.25,.25) node {\tiny $5$};
\draw (.6,.15) node {\tiny $6$};
\filldraw [white]
(-2,0) circle (1.5pt)
(0,0) circle (1.5pt)
(2,0) circle (1.5pt)
(-2,2) circle (1.5pt)
(0,2) circle (1.5pt)
(2,2) circle (1.5pt);
\draw[blue] (-1.5,0)--(-1.5,2);
\draw[blue] (-1.4,.5) node{\tiny $a$};
\end{tikzpicture}
\ar[r]^{T_{\check23}}&
\begin{tikzpicture}[
baseline=(x.base)
]
\node (x) at (0,.8){};
\draw[fill=blue!10] (-2,0) rectangle (2,2);
\draw[fill=white] (0,2) .. controls (-1.5,1.5) and (-0.5,0.5) .. (0,2);
\draw (0,2) .. controls (-.9,1.7) and (-1.7,.9) .. (-2,0);
\draw (0,2) .. controls (-.9,1.7) and (-1.6,1) .. (0,0);
\draw[fill=white] (2,2) .. controls (.5,1.5) and (1.5,0.5) .. (2,2);
\draw (2,2) .. controls (1.1,1.7) and (.3,.9) .. (0,0);
\draw (2,2) .. controls (1.7,1.1) and (.9,.3) .. (0,0);
\draw (0,2)--(0,0);
\draw (-1.9,.5) node {\tiny $1$};
\draw (-.1,.25) node {\tiny $2$};
\draw (-1.8,.15) node {\tiny $3$};
\draw (.1,.5) node {\tiny $4$};
\draw (.25,.25) node {\tiny $5$};
\draw (.6,.15) node {\tiny $6$};
\filldraw [white]
(-2,0) circle (1.5pt)
(0,0) circle (1.5pt)
(2,0) circle (1.5pt)
(-2,2) circle (1.5pt)
(0,2) circle (1.5pt)
(2,2) circle (1.5pt);
\draw[blue] (-1.5,0)--(-1.5,2);
\draw[blue] (-1.4,.5) node{\tiny $a$};
\end{tikzpicture}
}
\]
where the operator above the arrow realizes the corresponding element of the groupoid of decorated ideal triangulations within the quantum Teichm\"uller theory.
Thus,
\[
\zeta^{-6}\bar\fun_a=\Ad(T_{\check23})(T_{1\hat3})=T_{\check23}T_{1\hat3}\bar T_{\check23}=
T_{1\hat3}T_{1\hat2},
\]
where in the last equality, we have applied once the Pentagon relation, and we use the notation $\bar T=T^{-1}$. Here, we use the normalization where the braid-type and the lantern relations are satisfied without projective factors. By the above mentioned left-right symmetry $(1,2,3)\leftrightarrow (4,5,6)$, we immediately get the quantum realization of the Dehn twist $D_c$:
\[
\zeta^{-6}\bar\fun_c=T_{4\hat6}T_{4\hat5}.
\]
To calculate the quantum realization of $D_b$ we use a two-step chain of transformations of $\tau$:
\[
\xymatrix{
\begin{tikzpicture}[
]
\node (x) at (0,.8){};
\draw[fill=blue!10] (-2,0) rectangle (2,2);
\draw[fill=white] (0,2) .. controls (-1.5,1.5) and (-0.5,0.5) .. (0,2);
\draw[fill=white] (0,2) .. controls (-1.5,1.5) and (-0.5,0.5) .. (0,2);
\draw (0,2) .. controls (-.9,1.7) and (-1.7,.9) .. (-2,0);
\draw (0,2) .. controls (-.3,1.1) and (-1.1,.3) .. (-2,0);
\draw[fill=white] (2,2) .. controls (.5,1.5) and (1.5,0.5) .. (2,2);
\draw (2,2) .. controls (1.1,1.7) and (.3,.9) .. (0,0);
\draw (2,2) .. controls (1.7,1.1) and (.9,.3) .. (0,0);
\draw (-1.9,.5) node {\tiny $1$};
\draw (-1.75,.25) node {\tiny $2$};
\draw (-1.4,.15) node {\tiny $3$};
\draw (.1,.5) node {\tiny $4$};
\draw (.25,.25) node {\tiny $5$};
\draw (.6,.15) node {\tiny $6$};
\draw (0,2)--(0,0);
\filldraw [white]
(-2,0) circle (1.5pt)
(0,0) circle (1.5pt)
(2,0) circle (1.5pt)
(-2,2) circle (1.5pt)
(0,2) circle (1.5pt)
(2,2) circle (1.5pt);
\draw[blue] (-1,0)--(1,2);
\draw[blue] (1,0)--(2,1);
\draw[blue] (-2,1)--(-1,2);
\draw[blue] (-.3,.5) node{\tiny $b$};
\end{tikzpicture}
\ar[r]^{\bar T_{64}}
&
\begin{tikzpicture}[
]
\node (x) at (0,.8){};
\fill[blue!10] (-2,0) rectangle (2,2);
\draw[fill=white] (0,2) .. controls (-1.5,1.5) and (-0.5,0.5) .. (0,2);
\draw[fill=white] (0,2) .. controls (-1.5,1.5) and (-0.5,0.5) .. (0,2);
\draw (0,2) .. controls (-.9,1.7) and (-1.7,.9) .. (-2,0);
\draw (0,2) .. controls (-.3,1.1) and (-1.1,.3) .. (-2,0);
\draw[fill=white] (2,2) .. controls (.5,1.5) and (1.5,0.5) .. (2,2);
\draw (2,2) .. controls (1.1,1.7) and (.3,.9) .. (0,0);
\draw (2,2) .. controls (1.7,1.1) and (.9,.3) .. (0,0);
\draw (-1.9,.5) node {\tiny $1$};
\draw (-1.75,.25) node {\tiny $2$};
\draw (-1.4,.15) node {\tiny $3$};
\draw (.68,1.3) node {\tiny $4$};
\draw (.25,.25) node {\tiny $5$};
\draw (.1,1.85) node {\tiny $6$};
\draw (0,2)--(0,0);
\draw (-2,2)--(-2,0);
\draw (2,2)--(2,0);
\draw (-2,2)--(0,2);
\draw (-2,0)--(0,0);
\draw (0,0) .. controls (.3,.9) and (1,1.9) .. (1,2);
\draw (1,0) .. controls (1,.1) and (1.7,1.1) .. (2,2);
\filldraw [white]
(-2,0) circle (1.5pt)
(0,0) circle (1.5pt)
(2,0) circle (1.5pt)
(-2,2) circle (1.5pt)
(0,2) circle (1.5pt)
(2,2) circle (1.5pt);
\draw[blue] (-1,0)--(1,2);
\draw[blue] (1,0)--(2,1);
\draw[blue] (-2,1)--(-1,2);
\draw[blue] (-.3,.5) node{\tiny $b$};
\end{tikzpicture}
\ar[r]^{T_{41}\bar T_{63}}
&
\begin{tikzpicture}[
]
\node (x) at (0,.8){};
\fill[blue!10] (-2,0) rectangle (2,2);
\draw[fill=white] (0,2) .. controls (-1.5,1.5) and (-0.5,0.5) .. (0,2);
\draw[fill=white] (0,2) .. controls (-1.5,1.5) and (-0.5,0.5) .. (0,2);
\draw (0,2) .. controls (-.9,1.7) and (-1.7,.9) .. (-2,0);
\draw (0,2) .. controls (-.3,1.1) and (-1.1,.3) .. (-2,0);
\draw[fill=white] (2,2) .. controls (.5,1.5) and (1.5,0.5) .. (2,2);
\draw (2,2) .. controls (1.1,1.7) and (.3,.9) .. (0,0);
\draw (2,2) .. controls (1.7,1.1) and (.9,.3) .. (0,0);
\draw (1.9,.15) node {\tiny $1$};
\draw (-1.75,.25) node {\tiny $2$};
\draw (-1.4,.12) node {\tiny $3$};
\draw (.8,1.5) node {\tiny $4$};
\draw (.25,.25) node {\tiny $5$};
\draw (.1,1.85) node {\tiny $6$};
\draw (-2,2)--(0,2);
\draw (-2,0)--(0,0);
\draw (0,0) .. controls (.3,.9) and (1,1.9) .. (1,2);
\draw (1,0) .. controls (1,.1) and (1.7,1.1) .. (2,2);
\draw (0,0) .. controls (.4,1.2) and (1.4,1.9) .. (1.5,2);
\draw (1.5,0)--(2,.5);
\draw (0,2) .. controls (-1.2,1.6) and (-1.8,.7) .. (-2,.5);
\draw (-2,0) .. controls (-.5,.5) and (.4,1.9) .. (.5,2);
\draw (.5,0) .. controls (.6,.1) and (1.6,.8) .. (2,2);
\filldraw [white]
(-2,0) circle (1.5pt)
(0,0) circle (1.5pt)
(2,0) circle (1.5pt)
(-2,2) circle (1.5pt)
(0,2) circle (1.5pt)
(2,2) circle (1.5pt);
\draw[blue] (-1,0)--(1,2);
\draw[blue] (1,0)--(2,1);
\draw[blue] (-2,1)--(-1,2);
\draw[blue] (-.3,.5) node{\tiny $b$};
\end{tikzpicture}
}
\]
Thus, we have the following sequence of equalities:
\begin{multline*}
\zeta^{-6}\bar\fun_b= \Ad(\bar T_{64}T_{41}\bar T_{63})(T_{34})=
\bar T_{64}T_{41}\underline{\bar T_{63}T_{34}T_{63}}\bar T_{41}T_{64}\\=
\underline{\bar T_{64}T_{41}T_{64}}T_{34}\bar T_{41}T_{64}
=T_{61}\underline{T_{41}T_{34}\bar T_{41}}T_{64}=T_{61}T_{34}T_{31}T_{64},
\end{multline*}
where in each step the underlined fragment is transformed by using the Pentagon relation.

\vspace{0.2cm}
\noindent
To calculate the realization of $D_e$, we consider the following sequence of ideal triangulations:
\[
\xymatrix{
\begin{tikzpicture}[
baseline=(x.base),
]
\node (x) at (0,.8){};
\fill[blue!10] (-2,0) rectangle (2,2);
\draw[fill=white] (0,2) .. controls (-1.5,1.5) and (-0.5,0.5) .. (0,2);
\draw (0,2) .. controls (-.9,1.7) and (-1.7,.9) .. (-2,0);
\draw (0,2) .. controls (-.3,1.1) and (-1.1,.3) .. (-2,0);
\draw[fill=white] (2,2) .. controls (.5,1.5) and (1.5,0.5) .. (2,2);
\draw (2,2) .. controls (1.1,1.7) and (.3,.9) .. (0,0);
\draw (2,2) .. controls (1.7,1.1) and (.9,.3) .. (0,0);
\draw (0,2)--(0,0);
\draw (-2,2)--(-2,0);
\draw (2,2)--(2,0);
\draw (-2,2)--(0,2);
\draw (-2,0)--(0,0);
\draw (0,0)--(2,0);
\draw (0,2)--(2,2);
\draw (-1.9,.5) node {\tiny $1$};
\draw (-1.75,.25) node {\tiny $2$};
\draw (-1.4,.15) node {\tiny $3$};
\draw (.1,.5) node {\tiny $4$};
\draw (.25,.25) node {\tiny $5$};
\draw (.6,.15) node {\tiny $6$};
\filldraw [white]
(-2,0) circle (1.5pt)
(0,0) circle (1.5pt)
(2,0) circle (1.5pt)
(-2,2) circle (1.5pt)
(0,2) circle (1.5pt)
(2,2) circle (1.5pt);
\draw[blue] (0.2,0) arc (0:90:.2);
\draw[blue] (0,0.2) arc (90:180:.5 and .2);
\draw[blue] (0.2,2) arc (0:-90:.2 and .5);
\draw[blue] (-.5,2) .. controls (-.5,1.8) and (-1.4,1.4).. (-1,1);
\draw[blue] (0,1.5) .. controls (-.2,1.5) and (-.6,.6).. (-1,1);
\draw[blue] (-.9,1.1) node{\tiny $e$};
\end{tikzpicture}\ar[d]^{T_{34}}
&&
\begin{tikzpicture}[
baseline=(x.base),
]
\node (x) at (0,.8){};
\fill[blue!10] (-2,0) rectangle (2,2);
\draw[fill=white] (0,2) .. controls (-1.5,1.5) and (-0.5,0.5) .. (0,2);
\draw (0,2) .. controls (-.9,1.7) and (-1.7,.9) .. (-2,0);
\draw (0,2) .. controls (-.3,1.1) and (-1.1,.3) .. (-2,0);
\draw[fill=white] (2,2) .. controls (.5,1.5) and (1.5,0.5) .. (2,2);
\draw (-2,2)--(-2,0);
\draw (2,2)--(2,0);
\draw (-2,0)..controls(-1.5,2) and (-1.7,2)..(-1.5,2);
\draw (-1.5,0)..controls(-.72,0) and (.8,1.6)..(2,2);
\draw (-2,0)..controls(-.5,0) and (.5,1.5)..(2,2);
\draw (-2,0)..controls(-.8,.4) and (.8,1.8)..(1,2);
\draw (1,0)..controls(1.2,.2) and (2,1)..(2,2);
\draw (-2,0)..controls(-1.6,1.6) and (-1.1,1.9)..(-1,2);
\draw (-1,0)..controls(0,1) and (1.43,.29)..(2,2);
\draw (-2,0)..controls(-.8,.4) and (.5,1.8)..(.5,2);
\draw (-2,0)..controls(-1.6,1.4) and (-.7,1.9)..(-.7,2);
\draw (-.7,0)..controls(-.7,.5) and (.5,.5)..(.5,0);
\draw (-1.92,.75) node {\tiny $1$};
\draw (-1.75,.25) node {\tiny $2$};
\draw (-1.12,.29) node {\tiny $3$};
\draw (1.22,1.12) node {\tiny $4$};
\draw (.36,1.58) node {\tiny $5$};
\draw (0.1,1.87) node {\tiny $6$};
\filldraw [white]
(-2,0) circle (1.5pt)
(0,0) circle (1.5pt)
(2,0) circle (1.5pt)
(-2,2) circle (1.5pt)
(0,2) circle (1.5pt)
(2,2) circle (1.5pt);
\draw[blue] (0.2,0) arc (0:90:.2);
\draw[blue] (0,0.2) arc (90:180:.5 and .2);
\draw[blue] (0.2,2) arc (0:-90:.2 and .5);
\draw[blue] (-.5,2) .. controls (-.5,1.8) and (-1.4,1.4).. (-1,1);
\draw[blue] (0,1.5) .. controls (-.2,1.5) and (-.6,.6).. (-1,1);
\draw[blue] (-.9,1.1) node{\tiny $e$};
\end{tikzpicture}
\\
\begin{tikzpicture}[
baseline=(x.base),
]
\node (x) at (0,.8){};
\fill[blue!10] (-2,0) rectangle (2,2);
\draw[fill=white] (0,2) .. controls (-1.5,1.5) and (-0.5,0.5) .. (0,2);
\draw (0,2) .. controls (-.9,1.7) and (-1.7,.9) .. (-2,0);
\draw (0,2) .. controls (-.3,1.1) and (-1.1,.3) .. (-2,0);
\draw[fill=white] (2,2) .. controls (.5,1.5) and (1.5,0.5) .. (2,2);
\draw (2,2) .. controls (1.1,1.7) and (.3,.9) .. (0,0);
\draw (2,2) .. controls (1.7,1.1) and (.9,.3) .. (0,0);
\draw (-2,2)--(-2,0);
\draw (2,2)--(2,0);
\draw (-2,2)--(0,2);
\draw (-2,0)--(0,0);
\draw (0,0)--(2,0);
\draw (0,2)--(2,2);
\draw (-2,0)..controls(-.5,0) and (.5,1.5)..(2,2);
\draw (-1.9,.5) node {\tiny $1$};
\draw (-1.75,.25) node {\tiny $2$};
\draw (-1.22,.26) node {\tiny $3$};
\draw (0,.32) node {\tiny $4$};
\draw (.25,.25) node {\tiny $5$};
\draw (.6,.15) node {\tiny $6$};
\filldraw [white]
(-2,0) circle (1.5pt)
(0,0) circle (1.5pt)
(2,0) circle (1.5pt)
(-2,2) circle (1.5pt)
(0,2) circle (1.5pt)
(2,2) circle (1.5pt);
\draw[blue] (0.2,0) arc (0:90:.2);
\draw[blue] (0,0.2) arc (90:180:.5 and .2);
\draw[blue] (0.2,2) arc (0:-90:.2 and .5);
\draw[blue] (-.5,2) .. controls (-.5,1.8) and (-1.4,1.4).. (-1,1);
\draw[blue] (0,1.5) .. controls (-.2,1.5) and (-.6,.6).. (-1,1);
\draw[blue] (-.9,1.1) node{\tiny $e$};
\end{tikzpicture}\ar[r]^{T_{14}\bar T_{63}}
&
\begin{tikzpicture}[
baseline=(x.base),
]
\node (x) at (0,.8){};
\fill[blue!10] (-2,0) rectangle (2,2);
\draw[fill=white] (0,2) .. controls (-1.5,1.5) and (-0.5,0.5) .. (0,2);
\draw (0,2) .. controls (-.9,1.7) and (-1.7,.9) .. (-2,0);
\draw (0,2) .. controls (-.3,1.1) and (-1.1,.3) .. (-2,0);
\draw[fill=white] (2,2) .. controls (.5,1.5) and (1.5,0.5) .. (2,2);
\draw (2,2) .. controls (1.1,1.7) and (.3,.9) .. (0,0);
\draw (2,2) .. controls (1.7,1.1) and (.9,.3) .. (0,0);
\draw (-2,2)--(-2,0);
\draw (2,2)--(2,0);
\draw (-2,0)..controls(-1.5,2) and (-1.7,2)..(-1.5,2);
\draw (-1.5,0)..controls(-.72,0) and (.8,1.6)..(2,2);
\draw (-2,0)..controls(-.5,0) and (.5,1.5)..(2,2);
\draw (-2,0)..controls(-.8,.4) and (.8,1.8)..(1,2);
\draw (1,0)..controls(1.2,.2) and (2,1)..(2,2);
\draw (-1.92,.75) node {\tiny $1$};
\draw (-1.75,.25) node {\tiny $2$};
\draw (-1.12,.29) node {\tiny $3$};
\draw (0,.32) node {\tiny $4$};
\draw (.25,.25) node {\tiny $5$};
\draw (.6,.15) node {\tiny $6$};
\filldraw [white]
(-2,0) circle (1.5pt)
(0,0) circle (1.5pt)
(2,0) circle (1.5pt)
(-2,2) circle (1.5pt)
(0,2) circle (1.5pt)
(2,2) circle (1.5pt);
\draw[blue] (0.2,0) arc (0:90:.2);
\draw[blue] (0,0.2) arc (90:180:.5 and .2);
\draw[blue] (0.2,2) arc (0:-90:.2 and .5);
\draw[blue] (-.5,2) .. controls (-.5,1.8) and (-1.4,1.4).. (-1,1);
\draw[blue] (0,1.5) .. controls (-.2,1.5) and (-.6,.6).. (-1,1);
\draw[blue] (-.9,1.1) node{\tiny $e$};
\end{tikzpicture}\ar[r]^{T_{\check45}}
&
\begin{tikzpicture}[
baseline=(x.base),
]
\node (x) at (0,.8){};
\fill[blue!10] (-2,0) rectangle (2,2);
\draw[fill=white] (0,2) .. controls (-1.5,1.5) and (-0.5,0.5) .. (0,2);
\draw (0,2) .. controls (-.9,1.7) and (-1.7,.9) .. (-2,0);
\draw (0,2) .. controls (-.3,1.1) and (-1.1,.3) .. (-2,0);
\draw[fill=white] (2,2) .. controls (.5,1.5) and (1.5,0.5) .. (2,2);
\draw (2,2) .. controls (1.7,1.1) and (.9,.3) .. (0,0);
\draw (-2,2)--(-2,0);
\draw (2,2)--(2,0);
\draw (-2,0)..controls(-1.5,2) and (-1.7,2)..(-1.5,2);
\draw (-1.5,0)..controls(-.72,0) and (.8,1.6)..(2,2);
\draw (-2,0)..controls(-.5,0) and (.5,1.5)..(2,2);
\draw (-2,0)..controls(-.8,.4) and (.8,1.8)..(1,2);
\draw (1,0)..controls(1.2,.2) and (2,1)..(2,2);
\draw (-2,0)..controls(-1.6,1.6) and (-1.1,1.9)..(-1,2);
\draw (-1,0)..controls(0,1) and (1.43,.29)..(2,2);
\draw (-1.92,.75) node {\tiny $1$};
\draw (-1.75,.25) node {\tiny $2$};
\draw (-1.12,.29) node {\tiny $3$};
\draw (1.22,1.12) node {\tiny $4$};
\draw (.25,.25) node {\tiny $5$};
\draw (.6,.15) node {\tiny $6$};
\filldraw [white]
(-2,0) circle (1.5pt)
(0,0) circle (1.5pt)
(2,0) circle (1.5pt)
(-2,2) circle (1.5pt)
(0,2) circle (1.5pt)
(2,2) circle (1.5pt);
\draw[blue] (0.2,0) arc (0:90:.2);
\draw[blue] (0,0.2) arc (90:180:.5 and .2);
\draw[blue] (0.2,2) arc (0:-90:.2 and .5);
\draw[blue] (-.5,2) .. controls (-.5,1.8) and (-1.4,1.4).. (-1,1);
\draw[blue] (0,1.5) .. controls (-.2,1.5) and (-.6,.6).. (-1,1);
\draw[blue] (-.9,1.1) node{\tiny $e$};
\end{tikzpicture}
\ar[u]^{T_{\check56}}
}
\]
Thus, we have
\begin{multline*}
\zeta^{-6}\bar\fun_e= \Ad(T_{34}T_{14}\bar T_{63}T_{\check45}T_{\check56})(T_{\check2\hat6})=
T_{34}T_{14}\bar T_{63}T_{\check45}\underline{T_{\hat6\hat5}T_{\check2\hat6}\bar T_{\hat6\hat5}}\bar T_{\check45} T_{63}\bar T_{14}\bar T_{34}\\
=T_{34}T_{14}\bar T_{63}\underline{T_{\hat5\hat4}}T_{\check2\hat6}\underline{T_{\check2\hat5}\bar T_{\hat5\hat4}} T_{63}\bar T_{14}\bar T_{34}
=T_{34}\underline{T_{\hat4\check1}}\bar T_{63}T_{\check2\hat6}T_{\check2\hat5}\underline{T_{\check2\hat4}}T_{63}\underline{\bar T_{\hat4\check1}}\bar T_{34}\\
=T_{34}\underline{\bar T_{\hat3\check6}T_{\check6\hat2}}T_{\check2\hat5}T_{\check2\hat4}T_{\check2\check1}
\underline{T_{\hat3\check6}}\bar T_{34}=\underline{T_{34}T_{\check23}}T_{\check6\hat2}T_{\check2\hat5}
T_{\check2\hat4}T_{\check2\check1}
\bar T_{34}\\
=T_{\check23}T_{\check24}\underline{T_{\hat4\check3}}T_{\check6\hat2}
T_{\check2\hat5}\underline{T_{\check2\hat4}}T_{\check2\check1}
\underline{\bar T_{\hat4\check3}}
=T_{\check23}T_{\check24}T_{\check2\hat6}T_{\check2\hat5}T_{\check2\hat4} T_{\check2\check3}T_{\check2\check1},
\end{multline*}
where, as before, in each step the underlined fragment is transformed by applying the Pentagon relation.
We use throughout these computations the fact that $T_{ij}$ and $T_{kl}$ 
commute if $\{i,j\}\cap \{k,l\}=\emptyset$.
Again, using the symmetry $(1,2,3)\leftrightarrow(4,5,6)$, we also have
\[
\zeta^{-6}\bar\fun_f=T_{\check56}T_{\check51}T_{\check5\hat3}T_{\check5\hat2}T_{\check5\hat1} T_{\check5\check6}T_{\check5\check4}.
\]
In order to check the Chain relation, we first calculate the following product:
\[
\zeta^{-18}\bar\fun_c\bar\fun_b\bar\fun_a
=T_{4\hat6}T_{4\hat5}T_{61}T_{34}\underline{T_{31}}T_{64}\underline{T_{1\hat3}}T_{1\hat2}
=T_{4\hat6}T_{4\hat5}T_{61}T_{34}T_{64}\zeta P_{(31\hat3)}T_{1\hat2}=\zeta T_{4\hat6}T_{4\hat5}T_{61}T_{34}T_{64}T_{\hat3\hat2} P_{(31\hat3)},
\]
where we have applied the Inversion relation to the underlined fragment. Next, we calculate
\begin{multline*}
\zeta^{-36}(\bar\fun_c\bar\fun_b\bar\fun_a)^2
=\zeta^2 T_{4\hat6}T_{4\hat5}T_{61}T_{34}T_{64}T_{\hat3\hat2} \underline{ P_{(31\hat3)}}T_{4\hat6}T_{4\hat5}\underline{T_{61}T_{34}}T_{64}\underline{T_{\hat3\hat2} P_{(31\hat3)}}\\
=\zeta^2 T_{4\hat6}T_{4\hat5}T_{61}T_{34}\underline{T_{64}}T_{\hat3\hat2} \underline{T_{4\hat6}T_{4\hat5}T_{6\hat3}T_{14}T_{64}}T_{\hat1\hat2}P_{(3\hat3)(1\hat1)}\\
=\zeta^3 T_{4\hat6}T_{4\hat5}\underline{T_{61}T_{34}T_{\hat3\hat2}T_{\hat6\hat5}}T_{4\hat3}
T_{1\hat6}T_{4\hat6}T_{\hat1\hat2}P_{(64\hat6)}P_{(3\hat3)(1\hat1)}\\
=\zeta^3 \underline{T_{4\hat6}T_{4\hat5} T_{\hat6\hat5}}T_{\check51}\underline{T_{61}}T_{\hat3\hat2}T_{\check24}
\underline{T_{34}T_{4\hat3}
T_{1\hat6}}T_{4\hat6}T_{\hat1\hat2}P_{(64\hat6)}P_{(3\hat3)(1\hat1)}\\
=
\zeta^5 T_{\hat6\hat5} T_{4\hat6}T_{\check51}T_{\hat3\hat2}T_{\check24}\underline{P_{(61\hat6)}P_{(34\hat3)}
T_{4\hat6}T_{\hat1\hat2}P_{(64\hat6)}P_{(3\hat3)(1\hat1)}}
=
\zeta^5 T_{\hat6\hat5}T_{4\hat6}T_{\check51}T_{\check23}T_{\check24}
T_{\hat3\hat1}T_{\check6\hat2}P_{(1\check63\hat4\check1)},
\end{multline*}
where each equality is obtained by transforming the underlined fragment  by applying  the Pentagon relation (twice in the forth and once in the fifth equalities), the Inversion relation (once in the third and twice in the fifth equalities), and the extended symmetric group action (in the second, the third, and the sixth equalities).
Finally, taking the square of the obtained identity, we have
\begin{multline*}
\zeta^{-72}(\bar\fun_c\bar\fun_b\bar\fun_a)^4
=\zeta^{10} T_{\hat6\hat5}T_{4\hat6}T_{\check51}T_{\check23}T_{\check24}
T_{\hat3\hat1}T_{\check6\hat2}\underline{P_{(1\check63\hat4\check1)}
T_{\hat6\hat5}T_{4\hat6}T_{\check51}T_{\check23}T_{\check24}
T_{\hat3\hat1}T_{\check6\hat2}P_{(1\check63\hat4\check1)}}\\
=
\zeta^{10} T_{\check56}T_{4\hat6}T_{\check51}T_{\check23}T_{\check24}
\underline{T_{\hat3\hat1}}T_{\check6\hat2}\underline{T_{\check5\hat3}}
T_{\hat1\check3}T_{\check5\check6}T_{\check2\hat4}
T_{\check2\hat1}T_{\check46}
T_{3\hat2}P_{(13\check1)}P_{(46\check4)}
\\
=
\zeta^{10} T_{\check56}T_{4\hat6}T_{\check51}T_{\check23}T_{\check24}T_{\check6\hat2}T_{\check5\hat3}
T_{\check5\hat1}\underline{T_{\hat3\hat1}
T_{\hat1\check3}}T_{\check5\check6}T_{\check2\hat4}
\underline{T_{\check2\hat1}}T_{\check46}
\underline{T_{3\hat2}P_{(13\check1)}}P_{(46\check4)}
\\
=\zeta^{11} T_{\check56}\underline{T_{4\hat6}}T_{\check51}T_{\check23}\underline{T_{\check24}}
T_{\check6\hat2}T_{\check5\hat3}
T_{\check5\hat1}T_{\check5\check6}T_{\check2\hat4}T_{\check2\check3}
T_{\check46}T_{1\hat2}
P_{(46\check4)}\\
=\zeta^{11} T_{\check56}T_{\check51}T_{\check23}T_{\check24}
T_{\check2\hat6}\underline{T_{4\hat6}
T_{\check6\hat2}}T_{\check5\hat3}
T_{\check5\hat1}T_{\check5\check6}\underline{T_{\check2\hat4}}
T_{\check46}T_{\check2\check3}T_{1\hat2}
P_{(46\check4)}\\
=\zeta^{11} T_{\check56}T_{\check51}T_{\check23}T_{\check24}
T_{\check2\hat6}T_{\check5\hat3}
T_{\check5\hat1}T_{\check2\hat4}\underline{T_{4\hat6}T_{\check5\check6}}
T_{\check46}T_{\check2\check3}T_{1\hat2}
P_{(46\check4)}
\\=
\zeta^{11} T_{\check56}T_{\check51}T_{\check23}T_{\check24}
T_{\check2\hat6}T_{\check5\hat3}
T_{\check5\hat1}T_{\check2\hat4}T_{\check5\check6}T_{\check5\check4}\underline{T_{4\hat6}
T_{\check46}}T_{\check2\check3}T_{1\hat2}
\underline{P_{(46\check4)}}
\\=
\zeta^{12} T_{\check56}T_{\check51}\underline{T_{\check23}}T_{\check24}
\underline{T_{\check2\hat6}T_{\check5\hat3}}
T_{\check5\hat1}\underline{T_{\check2\hat4}T_{\check5\check6}T_{\check5\check4}}
T_{\check2\check3}T_{1\hat2}
\\=
\zeta^{12} T_{\check56}T_{\check51}T_{\check5\hat3}T_{\check5\hat2}T_{\check23}\underline{T_{\check24}}
T_{\check5\hat1}T_{\check5\check6}\underline{T_{\check5\hat2}}T_{\check2\hat6}
\underline{T_{\check5\check4}}
T_{\check5\hat2}T_{\check2\hat4}T_{\check2\check3}T_{1\hat2}
\\=
\zeta^{12} \underline{T_{\check56}T_{\check51}T_{\check5\hat3}T_{\check5\hat2}}T_{\check23}
\underline{T_{\check5\hat1}T_{\check5\check6}T_{\check5\check4}}T_{\check24}
T_{\check2\hat6}
T_{\check5\hat2}T_{\check2\hat4}T_{\check2\check3}T_{1\hat2}=\bar\fun_f\bar\fun_e,
\end{multline*}
where each equality, except for the last one, is obtained by transforming the underlined fragment  by applying  the Pentagon relation (one time in the third, the fifth, the sixth, the seventh, the tenth, and three times in the ninth equalities), the Inversion relation (in the forth and the eighth equalities), and the extended symmetric group action (in the second, the forth, and the eighth equalities), while in the last equality the underlined (respectively the non-underlined) fragment corresponds to the operator $\bar\fun_f$ (respectively $\bar\fun_e$)
\end{proof}

\begin{lemma}\label{puncturel}
Suppose that $s\geq 4$. Then the lift of each puncture relation is $\zeta^6$.
\end{lemma}
\begin{proof} 
Observe first that the central element 
$P_i$ which is the lift of the puncture relation at the puncture $p_i$ 
is independent of the particular subsurface $S^1_{0,3}$. If we consider another 
subsurface, there exists a homeomorphism $\varphi:S^s_{g,r}\to S^s_{g,r}$ 
fixing the puncture $p_i$ and sending it to the initial subsurface, because 
the boundary components are non-separating.    
The new  puncture relation is then conjugate of $P_i$ by $\widetilde{\varphi}$ 
and hence they  coincide, as they are elements of the center. 

\vspace{0.2cm}\noindent 
If $s\geq 4$ then there is an embedding 
$S^4_{0,3}\subset S^s_{g,r}$, such that each boundary component 
of $S^4_{0,3}$ has a puncture on it.  
Consider first the following decomposition  $\tau$ 
of the punctured pair of pants 
into triangles. The position of the label of each triangle 
indicates also the marked corner. 

\vspace{0.2cm}
\begin{center}
\includegraphics[scale=0.3]{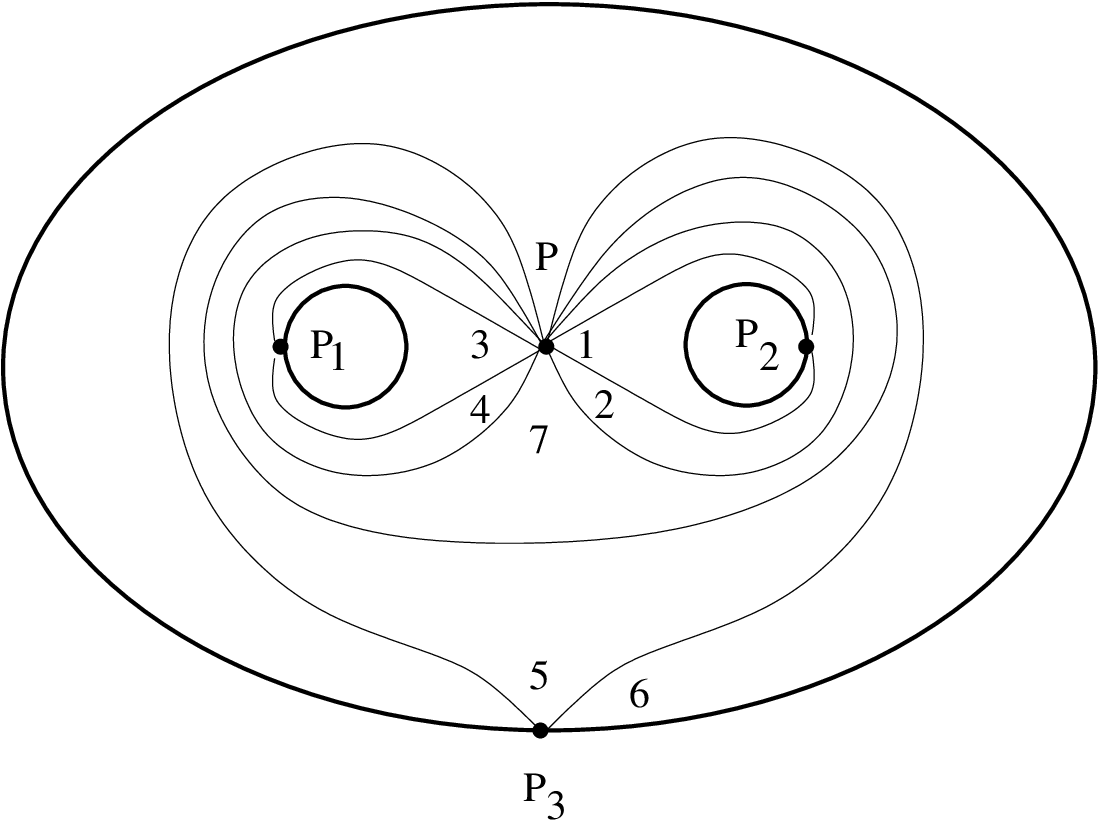}
\end{center}
\vspace{0.2cm}

\vspace{0.2cm}\noindent 
Then we can express easily the action of each Dehn twist $D_{a_j}$ 
on the triangulation $\tau$ as a composition of flips. 
If we set ${\mathsf F}_{a_j}={\mathsf F}(\tau, D_{a_j}(\tau))$ then we have:  

\[{\mathsf F}_{a_1}=T_{\check3\check4}^{-1}, \; {\mathsf F}_{a_2}=T_{\check12}^{-1},\;  
{\mathsf F}_{a_3}=T_{\check56}^{-1}\] 

\vspace{0.2cm}\noindent 
Further we use the sequence of transformations below, 
in order to change the triangulation $\tau$ into a triangulation 
which intersects the curve $a_{12}$ in only two points.

\vspace{0.2cm}
\begin{center}
\includegraphics[scale=0.26]{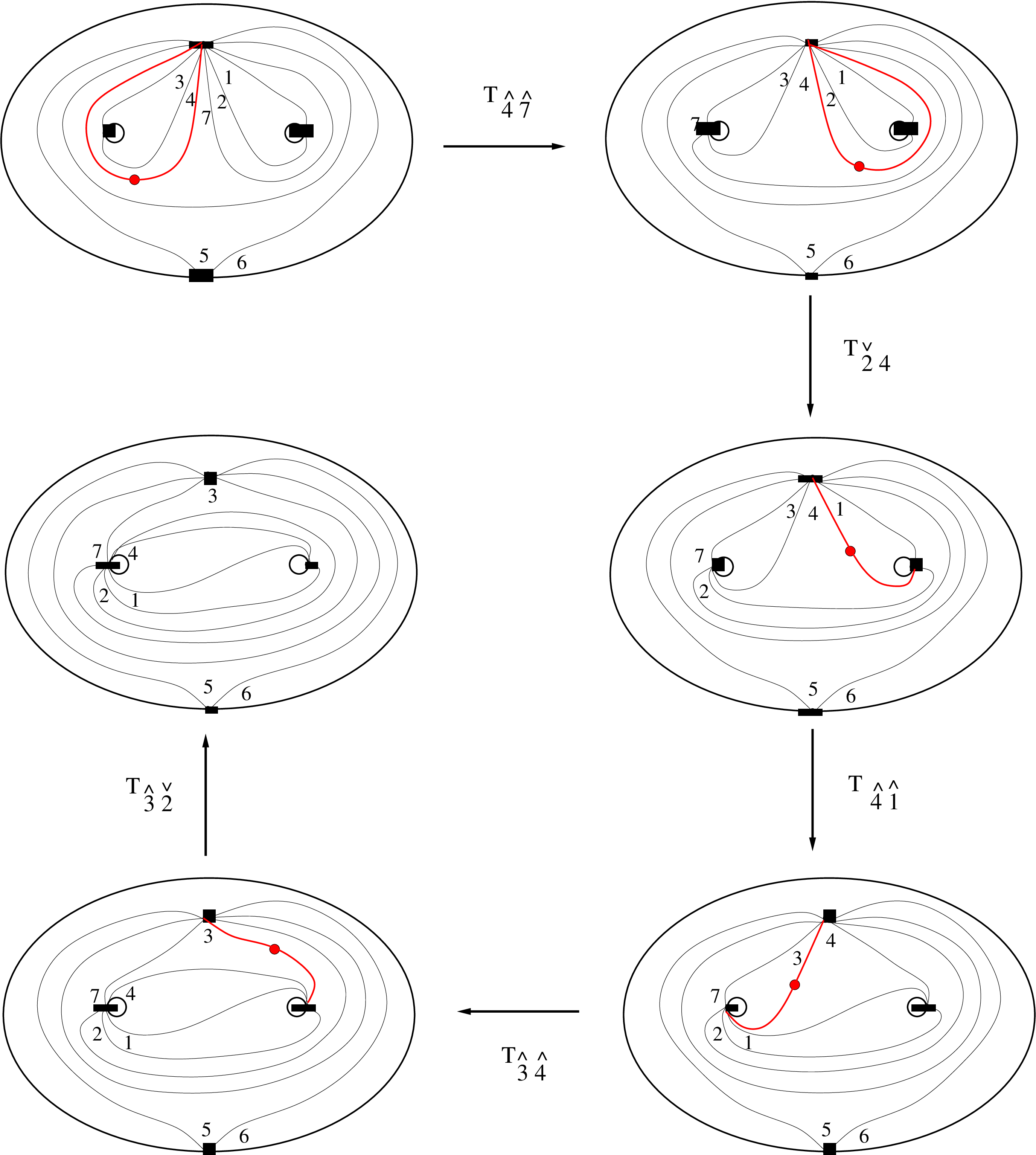}
\end{center}
\vspace{0.2cm}
\vspace{0.2cm}\noindent 
Here and in the pictures below we marked by a dot
the edges where a flip occurs, in order to help the reader
visualise the sequence of transformations.
Then the method outlined above permits to compute the Dehn twist 
${\mathsf F}_{a_{12}}={\mathsf F}(\tau, D_{a_{12}}(\tau))$ as follows: 

\[ {\mathsf F}_{a_{12}}= \Ad(T_{\hat{4}\hat{7}} T_{\check{2}4}T_{\hat{4}\hat{1}}T_{\hat{3}\hat{4}}T_{\hat{3}\check{2}})(T_{\check{3}\hat{7}}^{-1})\]

\vspace{0.2cm}\noindent 
Let us first simplify the formula for $\fun_{a_{12}}$. We have
\begin{multline*}
\bar\fun_{a_{12}}=T_{\check 74}T_{\check 24}T_{\hat4\hat1}T_{\hat3\hat4}
\underline{T_{\hat 3\check 2}T_{\check 7\hat 3}
\bar T_{\hat 3\check 2}}\bar T_{\hat3\hat4}\bar T_{\hat4\hat1}\bar T_{\check 24}\bar T_{\check 74}
=T_{\check 74}T_{\check 24}T_{\hat4\hat1}\underline{T_{\hat3\hat4}
T_{\check 7\hat 3}}T_{\check 7\check 2}\underline{\bar T_{\hat3\hat4}}\bar T_{\hat4\hat1}\bar T_{\check 24}\bar T_{\check 74}\\
=T_{\check 74}T_{\check 24}\underline{T_{\hat4\hat1}}T_{\check 7\hat 3}\underline{T_{\check 7\hat 4}}
T_{\check 7\check 2}\underline{\bar T_{\hat4\hat1}}\bar T_{\check 24}\bar T_{\check 74}
=T_{\check 74}T_{\check 24}T_{\check 7\hat 3}T_{\check 7\hat 4}T_{\check 7\hat 1}\underline{T_{\check 7\check 2}\bar T_{\check 24}\bar T_{\check 74}}\\
=T_{\check 74}\underline{T_{\check 24}}T_{\check 7\hat 3}\underline{T_{\check 7\hat 4}}T_{\check 7\hat 1}
\underline{\bar T_{\check 24}}T_{\check 7\check 2}
=T_{\check 74}T_{\check 7\hat 3}T_{\check 7\hat 4}T_{\check 7\hat2}T_{\check 7\hat 1}T_{\check 7\check 2}
\end{multline*}
where in each step the underlined fragment is transformed by using the Pentagon equation, and in the last equality it is also combined with the symmetry relation $T_{\check2 4}=T_{\hat 4\hat2}$.

\vspace{0.2cm}\noindent 
Our triangulation is invariant under the 
following simultaneous cyclic permutations
\[
\pi\colon P_1\mapsto P_2\mapsto P_3\mapsto P_1,\quad 1\mapsto\check6\mapsto3\mapsto 1,\quad 2\mapsto\hat5\mapsto\check4\mapsto2,\quad 7\mapsto\check 7,
\]
so that the contours $a_j$ and $a_{kl}$ are transformed as follows:
\[
\pi\colon a_1\mapsto a_2\mapsto a_3\mapsto a_1,\quad a_{12}\mapsto a_{23}\mapsto a_{31}\mapsto a_{12}.
\]
Thus, it suffices to know the explicit formula for $\fun_{a_{12}}$ in order to write out the other two without any further calculation:
\[
\bar\fun_{a_{23}}=\pi(\bar\fun_{a_{12}})=\pi(T_{\check 74}T_{\check 7\hat 3}T_{\check 7\hat 4}T_{\check 7\hat2}T_{\check 7\hat 1}T_{\check 7\check 2})=
T_{\hat 7\hat2}T_{\hat 7\hat 1}T_{\hat 7\check2}T_{\hat 7\check5}T_{\hat 76}T_{\hat 75},
\]
and 
\[
\bar\fun_{a_{31}}=\pi(\bar\fun_{a_{23}})=\pi(T_{\hat 7\hat2}T_{\hat 7\hat 1}T_{\hat 7\check2}T_{\hat 7\check5}T_{\hat 76}T_{\hat 75})=
T_{7\check5}T_{76}T_{75}T_{74}T_{7\hat3}T_{7\hat4}.
\]
Now, we have
\begin{multline*}
\bar\fun_{a_{12}}\bar\fun_{a_{23}}\bar\fun_{a_{31}}=
T_{\check 74}T_{\check 7\hat 3}T_{\check 7\hat 4}T_{\check 7\hat2}T_{\check 7\hat 1}\underline{T_{\check 7\check 2}T_{\hat 7\hat2}}T_{\hat 7\hat 1}T_{\hat 7\check2}T_{\hat 7\check5}T_{\hat 76}\underline{T_{\hat 75}
T_{7\check5}}T_{76}T_{75}T_{74}T_{7\hat3}T_{7\hat4}\\
=T_{\check 74}T_{\check 7\hat 3}T_{\check 7\hat 4}T_{\check 7\hat2}T_{\check 7\hat 1}\zeta P_{(2\hat 7\hat2)}T_{\hat 7\hat 1}T_{\hat 7\check2}T_{\hat 7\check5}T_{\hat 76}\zeta
P_{(\hat57\check5)}T_{76}T_{75}T_{74}T_{7\hat3}T_{7\hat4}\\
=\zeta^2 T_{\check 74}T_{\check 7\hat 3}T_{\check 7\hat 4}\underline{T_{\check 7\hat2}T_{\check 7\hat 1}T_{\hat 2\hat 1}}T_{\hat 27}\underline{T_{\hat 2\check5}T_{\hat 26}T_{\check56}}T_{\check5\check2}T_{\check54}T_{\check5\hat3}T_{\check5\hat4}P_{(2\hat 75\check2)}\\
=\zeta^2 T_{\check 74}T_{\check 7\hat 3}T_{\check 7\hat 4}T_{\hat 2\hat 1}\underline{T_{\check 7\hat2}T_{\hat 27}}T_{\check56}\underline{T_{\hat 2\check5}T_{\check5\check2}}T_{\check54}T_{\check5\hat3}T_{\check5\hat4}P_{(2\hat 75\check2)}\\
=\zeta^2 T_{\check 74}T_{\check 7\hat 3}T_{\check 7\hat 4}T_{\hat 2\hat 1}\zeta P_{(\check 7\hat 27)}T_{\check56}\zeta P_{(\hat 2\check5\check2)}T_{\check54}T_{\check5\hat3}T_{\check5\hat4} P_{(2\hat 75\check2)}
=\zeta^4 T_{\check 74}T_{\check 7\hat 3}\underline{T_{\check 7\hat 4}}T_{\hat 2\hat 1}T_{\check56}\underline{T_{\hat74}}T_{\hat7\hat3}T_{\hat7\hat4} P_{(7\hat7)}\\
=\zeta^4 T_{\check 74}T_{\check 7\hat 3}T_{\hat 2\hat 1}T_{\check56}\zeta P_{(\check4\hat74)}T_{\hat7\hat3}T_{\hat7\hat4} P_{(7\hat7)}
=\zeta^5 \underline{T_{\check 74}T_{\check 7\hat 3}}T_{\hat 2\hat 1}T_{\check56} \underline{T_{4\hat3}}T_{47} P_{(74\check7)}\\
=\zeta^5 T_{\hat 2\hat 1}T_{\check56}T_{4\hat3}\underline{T_{\check 74}T_{47}} P_{(74\check7)}
=\zeta^5 T_{\hat 2\hat 1}T_{\check56}T_{4\hat3}\zeta P_{(\check 747)} P_{(74\check7)}=\zeta^6 T_{\hat 2\hat 1}T_{\check56}T_{4\hat3}=\zeta^6\bar\fun_{a_{2}}\bar\fun_{a_{3}}\bar\fun_{a_{1}}
\end{multline*}
where in the underlined fragments the Pentagon equation is used twice in the forth and once in the ninth equalities, the Inversion relation is used twice in the second and the fifth, and once in the seventh and the tenth equalities, while in the third, sixth, eighth, and eleventh equalities the permutation operators are moved to the right and the powers of  $\zeta$, to the left.

\end{proof}

\begin{proof}[Proof of Lemma \ref{normalization}] 
The idea of the proof is to calculate the lift of the 
lantern relation. 
Consider the following decorated 
triangulation $\tau$ of the 4-holed disk with 4 punctures: 

\vspace{0.2cm}
\begin{center}
\includegraphics[scale=0.3]{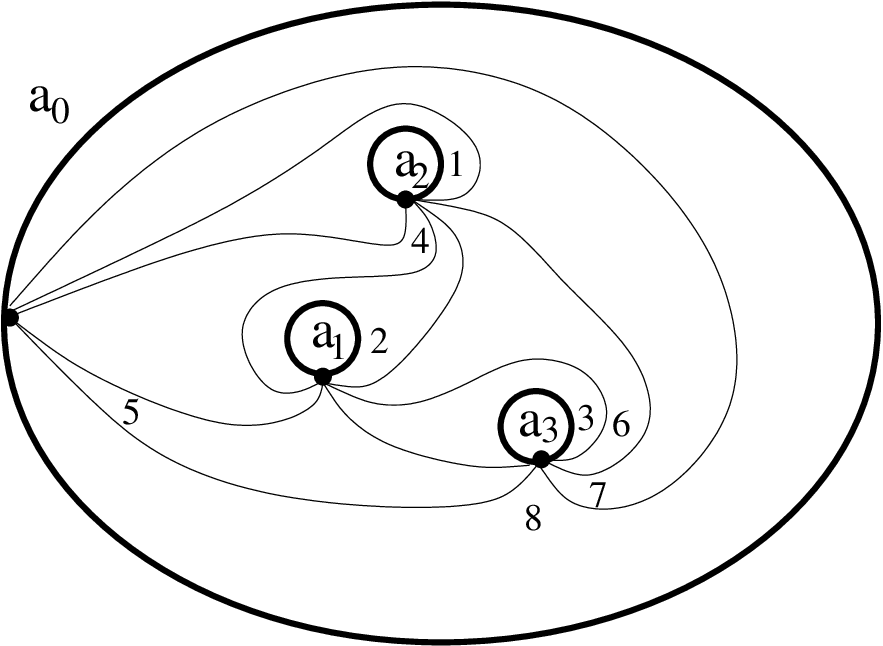}
\end{center}
The trick used in \cite{Ka,Ka2} for computing $D_{a}$ is to use 
a sequence of flips to change the triangulation into one which  
intersects some curve isotopic to $a$ into two points. Then 
the Dehn twist along $a$ can be expressed as the flip of one of the two 
edges of the latter triangulation  intersecting $a$. 
This recipe generalizes to the case where the curve $a$ intersects 
several edges of the triangulation, if $a$ is a boundary component with one puncture on it. 
Specifically, let $e_1,\ldots, e_s$ be the edges  issued 
from the puncture, in counterclockwise order. Then the Dehn twist $D_a$ 
can be expressed as the result of composing the flips of 
$e_1, e_2,\ldots, e_{s-1}$.  
We illustrate this procedure with the case of  the left Dehn twist 
$D^{-1}_{a_3}$ on the triangulation $\tau$ above:

\vspace{0.2cm}
\begin{center}
\includegraphics[scale=0.3]{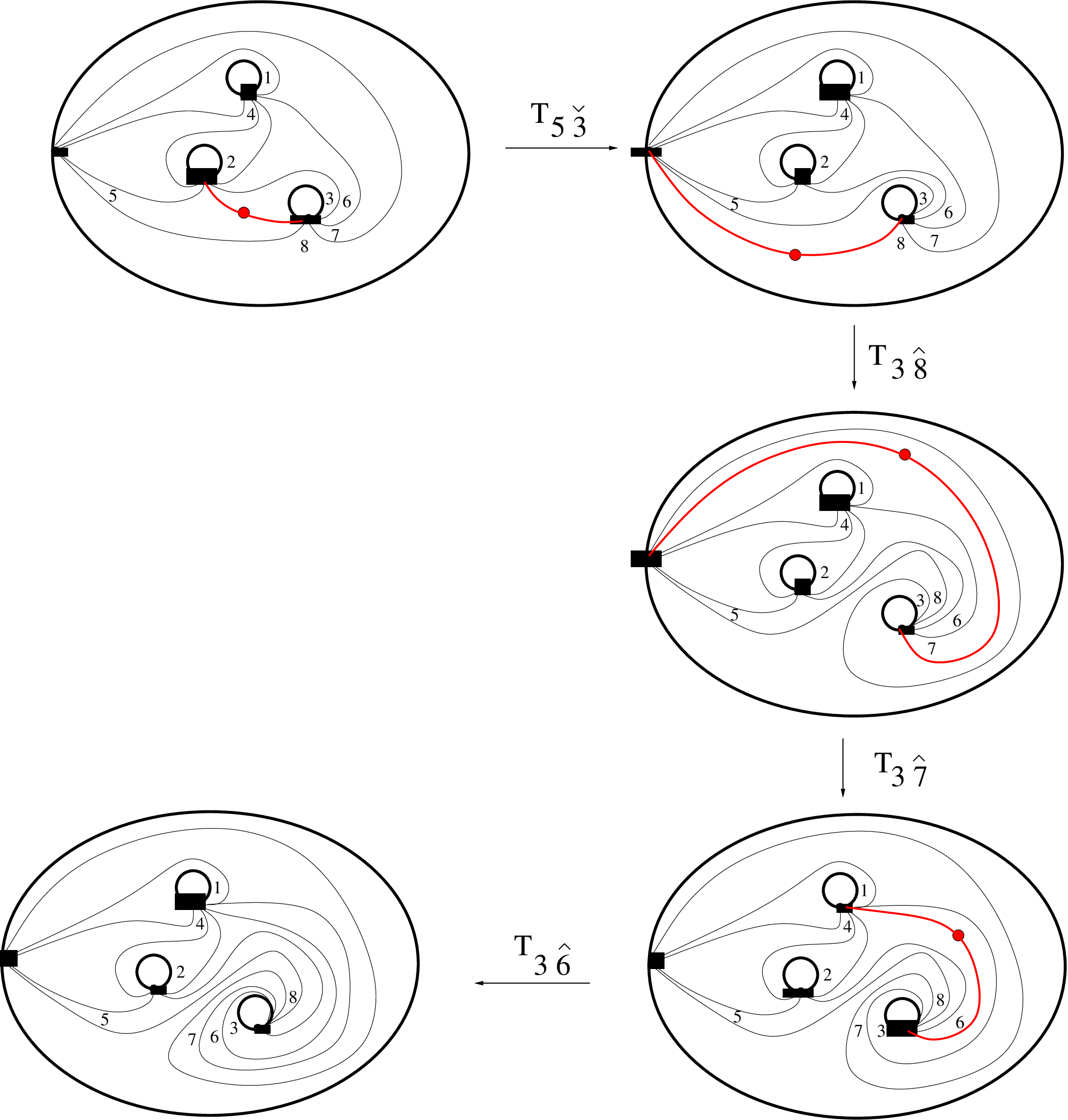}
\end{center}
In particular, we find the  following expression for 
the right Dehn twist along $a_3$: 
 \begin{equation}
\bar\fun_{a_3}=\bar \fun(\tau,D_{a_3}\tau)=T_{3\check5}T_{3\hat8}T_{3\hat7}T_{3\hat6}
\end{equation}
We used above the symmetry property of the $T$-matrix 
$T_{3\check5}=T_{5\check3}$ (see Remark \ref{symmetry} equation (\ref{symm})). 
The same recipe for the remaining Dehn twists along boundary 
components gives us:  
\begin{equation}
\bar\fun_{a_2}=\bar \fun(\tau,D_{a_2}\tau)=T_{24}T_{25}T_{2\check3}T_{26}
\end{equation}
 \begin{equation}
\bar\fun_{a_1}=\bar \fun(\tau,D_{a_1}\tau)=T_{1\hat4}T_{1\check2}T_{1\check6}T_{17}
\end{equation}
 \begin{equation}
\bar\fun_{a_0}=\bar \fun(\tau,D_{a_0}\tau)=T_{\check8\hat5}T_{\check8\check4}T_{\check8\check1}
T_{\check8\check7}
\end{equation}

\vspace{0.2cm}\noindent 
In order to compute $F_{a_{12}}$ we need to transform the triangulation 
$\tau$ into one which intersects a curve isotopic to $a_{12}$ into precisely
two points. 
This can be done as follows: 

\vspace{0.2cm}
\begin{center}
\includegraphics[scale=0.3]{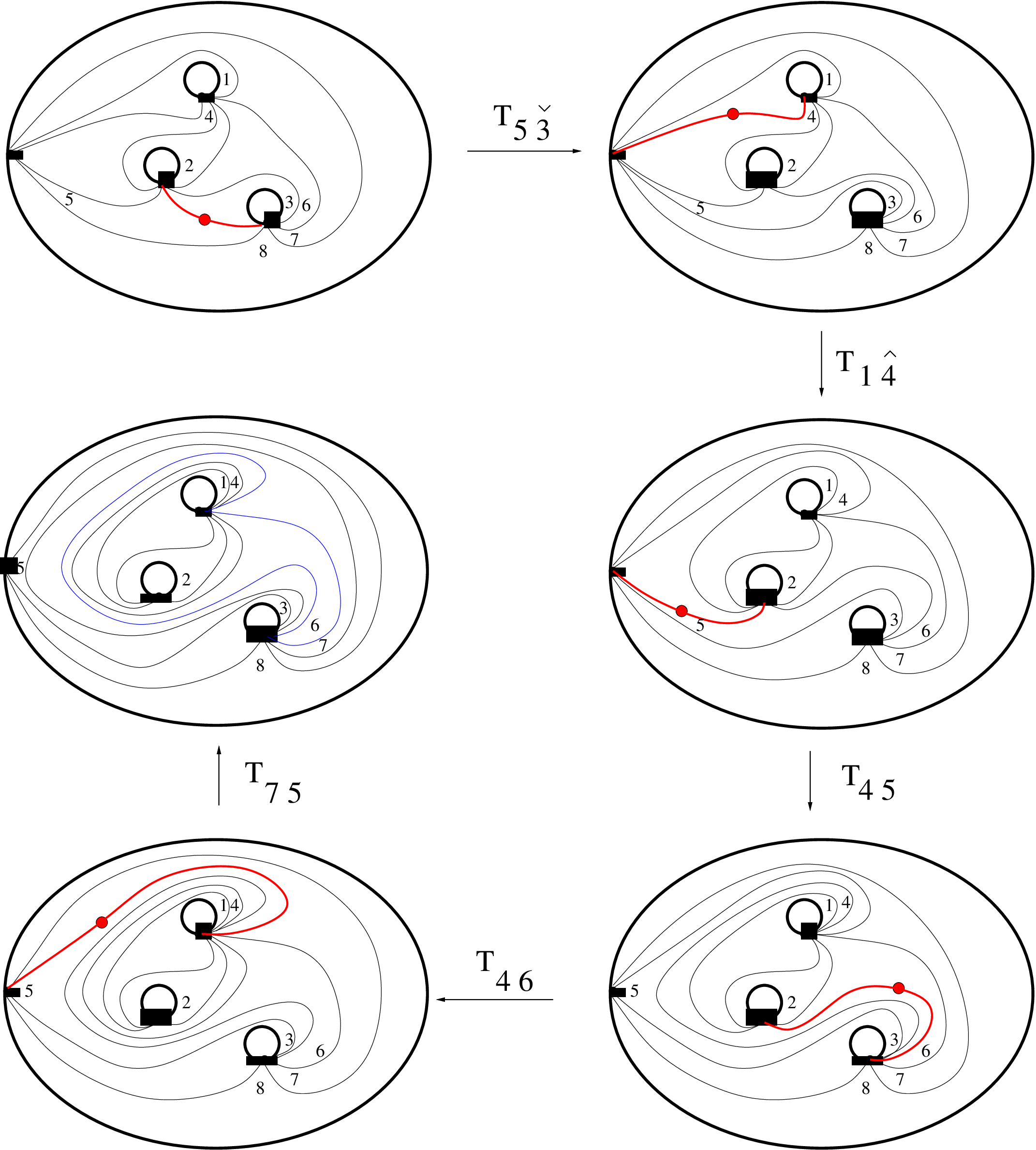}
\end{center}
Therefore we have: 
\begin{equation}
\bar\fun_{a_{12}}=\bar \fun(\tau,D_{a_{12}}\tau)=
\Ad(T_{3\check5}T_{1\check4}T_{45}T_{46}T_{75})(T_{\check67})
\end{equation}
The following  sequence of transformations 

\vspace{0.2cm}
\begin{center}
\includegraphics[scale=0.3]{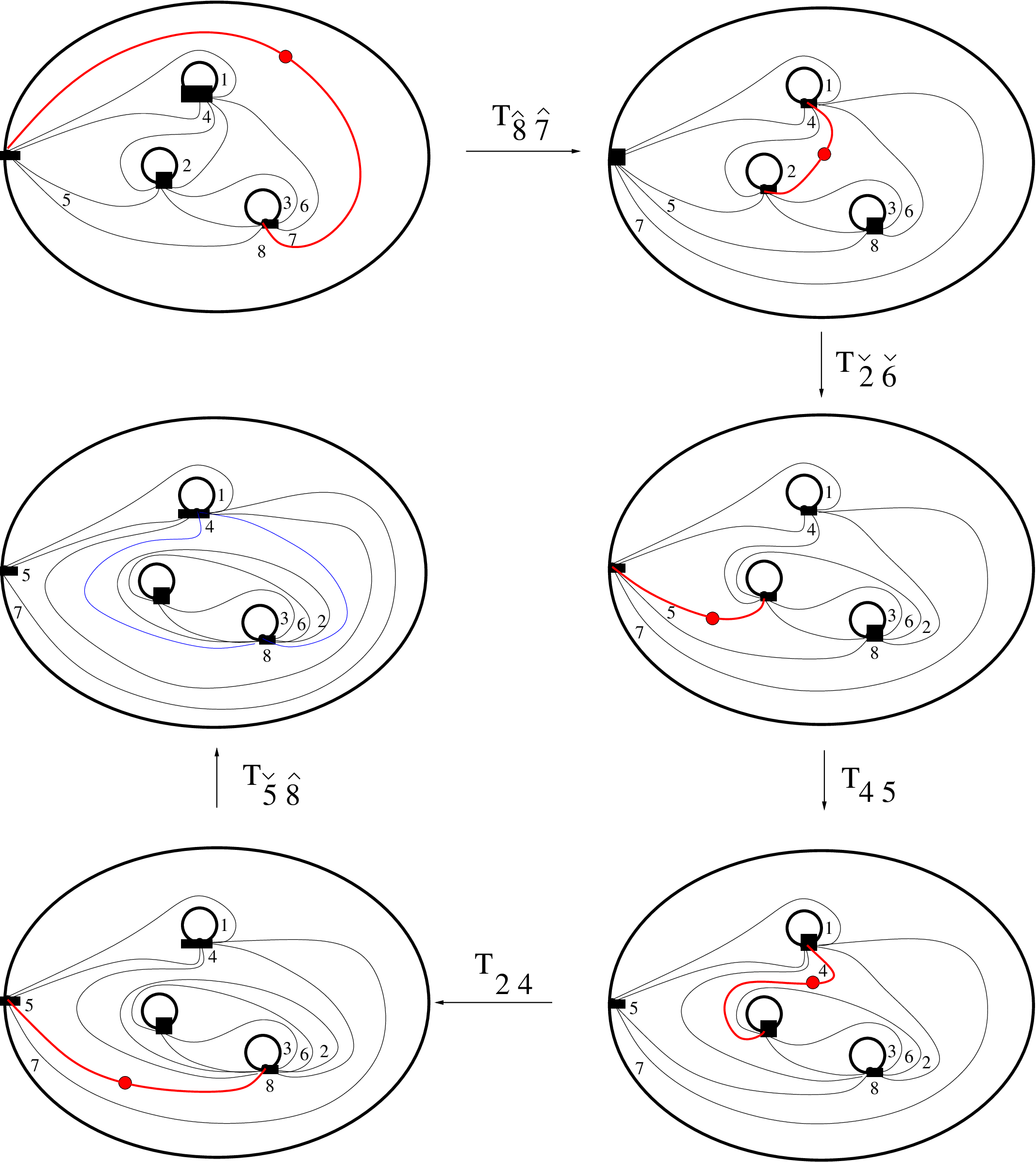}
\end{center}
can be used to compute:  
\begin{equation}
\bar\fun_{a_{13}}=\bar \fun(\tau,D_{a_{13}}\tau)=
\Ad(T_{\hat8\hat7}T_{\check2\check6}T_{45}T_{24}T_{\check5\hat8})(T_{\check4\check5})
\end{equation}
Eventually use the transformations  

\vspace{0.2cm}
\begin{center}
\includegraphics[scale=0.3]{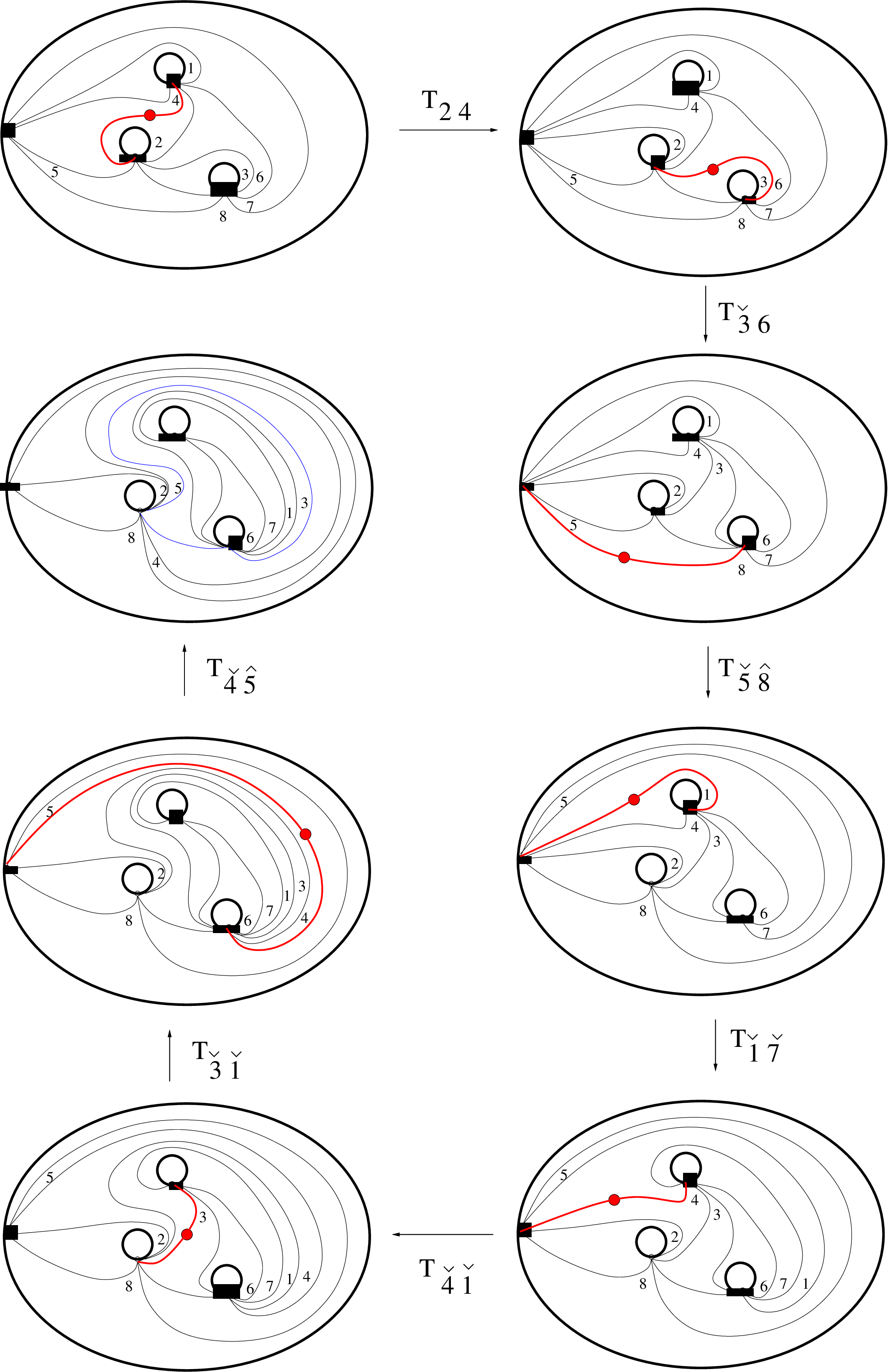}
\end{center}
in order to obtain: 
\begin{equation}
\bar\fun_{a_{23}}=\bar \fun(\tau,D_{a_{23}}\tau)=
\Ad(T_{24}T_{\check36}T_{\check5\hat8}T_{\check1\check7}T_{\check4\check1}T_{\check3\check1}T_{\check4\hat5})(T_{3\check5})
\end{equation}

\vspace{0.2cm}\noindent 
The next step is to simplify the expression of the last three Dehn twist, 
as follows: 
\begin{multline*}
\bar\fun_{a_{12}} =T_{3\check5}T_{1\check4}T_{45}\underline{T_{46}T_{75}}T_{\check67}
\underline{\bar T_{75}\bar T_{46}}\bar T_{45}\bar T_{1\check4} \bar T_{3\check5}= 
T_{3\check5}T_{1\check4}T_{45}T_{75}\underline{T_{46}T_{\check67}
\bar T_{46}}\bar T_{75}\bar T_{45}\bar T_{1\check4} \bar T_{3\check5}= \\
= T_{3\check5}T_{1\check4}T_{45}T_{75}T_{\check67}
\underline{\bar T_{47}\bar T_{75}\bar T_{45}}\bar T_{1\check4} \bar T_{3\check5}=
T_{3\check5}T_{1\check4}T_{45}\underline{T_{75}T_{\check67}
\bar T_{75}}T_{47}\bar T_{1\check4} \bar T_{3\check5}
= T_{3\check5}T_{1\check4}T_{45}T_{\check67}
T_{\check65}T_{47}\bar T_{1\check4} \bar T_{3\check5}=\\
=T_{1\check4}T_{\check67}T_{3\check5}T_{45}T_{\check65}\bar T_{3\check5}
T_{47}\bar T_{1\check4}=
T_{1\check4}T_{\check67}\underline{T_{3\check5}T_{45}\bar T_{3\check5}}\: 
\underline{T_{3\check5}T_{\check65}\bar T_{3\check5}}T_{47}\bar T_{1\check4}=
T_{1\check4}T_{\check67}T_{45}T_{3\check4}T_{\check65}T_{3\hat6}T_{47}\bar T_{1\check4}
\end{multline*}
The first equality above corresponds to the commutativity 
of $T_{ij}$ and $T_{kl}$ in the case when the 
two sets of indices are disjoint, for each one of the underlined fragments. 
We further also made use of the symmetry property 
from (Remark \ref{symmetry}, relation (\ref{symm}))
in order to be able to use the Pentagon relation, as 
in the last equality above. Specifically, the rightmost reduction 
consists of the the following steps:  
\begin{equation}
T_{3\check5}T_{\check65}\bar T_{3\check5}=
\underline{T_{5\check3}T_{\check65}}\bar T_{3\check5}=
\underline{T_{\check65}T_{3\hat6}T_{5\check3}}\bar T_{3\check5}=
T_{\check65}T_{3\hat6}T_{3\check5}\bar T_{3\check5}=
T_{\check65}T_{3\hat6}
\end{equation}

\vspace{0.2cm}\noindent 
Similar simplifications lead to:

\begin{multline*}
\bar\fun_{a_{13}} =T_{\hat8\hat7}T_{\check2\check6}T_{45}\underline{T_{24}T_{\check5\hat8}}
T_{\check4\check5}\underline{\bar T_{\check5\hat8}\bar T_{24}} \bar T_{45} \bar 
T_{\check2\check6}\bar T_{\hat8\hat7}=
T_{\hat8\hat7}T_{\check2\check6}T_{45}T_{\check5\hat8}\underline{T_{24}
T_{\check4\check5}\bar T_{24}}\bar T_{\check5\hat8} \bar T_{45} \bar 
T_{\check2\check6}\bar T_{\hat8\hat7}=\\
=T_{\hat8\hat7}T_{\check2\check6}\underline{T_{45}T_{\check5\hat8}
T_{\check4\check5}}T_{5\check2}\bar T_{\check5\hat8} \bar T_{45} \bar 
T_{\check2\check6}\bar T_{\hat8\hat7}= \zeta 
T_{\hat8\hat7}T_{\check2\check6}T_{\check5\hat8}
T_{4\hat8} P_{(45\hat4)}\bar T_{5\check2}\bar T_{\check5\hat8} \bar T_{45} \bar 
T_{\check2\check6}\bar T_{\hat8\hat7}=\\
=\zeta 
T_{\hat8\hat7}T_{\check2\check6}T_{\check5\hat8}
\underline{T_{4\hat8} \bar T_{\hat4\check2}\bar T_{4\hat8}} \bar T_{5\hat4} \bar 
T_{\check2\check6}\bar T_{\hat8\hat7}P_{(45\hat4)}=
\zeta 
T_{\hat8\hat7}\underline{T_{\check2\check6}T_{\check5\hat8}}
T_{24} T_{2\hat8}  \underline{\bar T_{5\hat4} \bar 
T_{\check2\check6}}\bar T_{\hat8\hat7}P_{(45\hat4)}=\\
=\zeta 
T_{\hat8\hat7}T_{\check5\hat8}T_{\check2\check6}
T_{24} T_{2\hat8} \bar 
T_{\check2\check6}\bar T_{5\hat4}  
\bar T_{\hat8\hat7}P_{(45\hat4)}=\zeta 
T_{\hat8\hat7}T_{\check5\hat8}\underline{T_{\check2\check6}
T_{24}\bar 
T_{\check2\check6}}\;\underline{T_{\check2\check6} T_{2\hat8} \bar 
T_{\check2\check6}}\bar T_{5\hat4}  
\bar T_{\hat8\hat7}P_{(45\hat4)}=\\
=\zeta 
T_{\hat8\hat7}T_{\check5\hat8}
T_{24}  
T_{\hat4\check6}
T_{2\hat8}  
T_{\check8\check6}
\bar T_{5\hat4}  
\bar T_{\hat8\hat7}P_{(45\hat4)}
\end{multline*}

\begin{multline*}
\bar\fun_{a_{23}} = 
T_{24}T_{\check 36}T_{\check5\hat8}T_{\check1\check7}\underline{T_{\check4\check1}T_{\check3\check1}}
T_{\check4\hat5}T_{3\check5}\underline{\bar T_{\check4\hat5}\bar T_{\check3\check1}}
\bar T_{\check4\check1}\bar T_{\check1\check7}\bar T_{\check5\hat8}\bar T_{\check 36}\bar T_{24}
=\\
=T_{24}T_{\check 36}T_{\check5\hat8}T_{\check1\check7}T_{\check4\check1}
T_{\check4\hat5}\underline{T_{\check3\check1}T_{3\check5}\bar T_{\check3\check1}}\bar T_{\check4\hat5}
\bar T_{\check4\check1}\bar T_{\check1\check7}\bar T_{\check5\hat8}\bar T_{\check 36}\bar T_{24}=\\
=
T_{24}T_{\check 36}T_{\check5\hat8}T_{\check1\check7}T_{\check4\check1}
T_{\check4\hat5}T_{3\check5}\underline{T_{5\check1}\bar T_{\check4\hat5}
\bar T_{\check4\check1}}\bar T_{\check1\check7}\bar T_{\check5\hat8}\bar T_{\check 36}\bar T_{24}= 
T_{24}T_{\check 36}T_{\check5\hat8}T_{\check1\check7}T_{\check4\check1}
\underline{T_{\check4\hat5}T_{3\check5}\bar T_{\check4\hat5}}
T_{5\check1}\bar T_{\check1\check7}\bar T_{\check5\hat8}\bar T_{\check 36}\bar T_{24}= 
\\
=T_{24}T_{\check 36}T_{\check5\hat8}T_{\check1\check7}T_{\check4\check1}
T_{\check4\hat5}T_{3\check5}\underline{T_{5\check1}\bar T_{\check4\hat5}
\bar T_{\check4\check1}}\bar T_{\check1\check7}\bar T_{\check5\hat8}\bar T_{\check 36}\bar T_{24}= 
T_{24}T_{\check 36}T_{\check5\hat8}T_{\check1\check7}T_{\check4\check1}
T_{3\check5}T_{3\hat4}
T_{5\check1}\bar T_{\check1\check7}\bar T_{\check5\hat8}\bar T_{\check 36}\bar T_{24}
\end{multline*}

\vspace{0.2cm}\noindent 
Putting all these together we obtain: 
\begin{multline*}
\bar\fun_{a_{12}}\bar\fun_{a_{23}}\bar\fun_{a_{13}}=
\zeta T_{1\check4}T_{\check67}T_{45}T_{3\check4}T_{\check65}T_{3\hat6}T_{47}\bar T_{1\check4}T_{24}T_{\check 36}T_{\check5\hat8}T_{\check1\check7}T_{\check4\check1}
T_{3\check5}T_{3\hat4}
T_{5\check1}\bar T_{\check1\check7}\underline{\bar T_{\check5\hat8}}\bar T_{\check 36}\underline{\bar T_{24}}
\underline{T_{\hat8\hat7}T_{\check5\hat8}}
\underline{T_{24}}  
T_{\hat4\check6}T_{2\hat8}T_{\check8\check6}\times \\  
\times   
\bar T_{5\hat4}\bar T_{\hat8\hat7}P_{(45\hat4)}=
\zeta T_{1\check4}T_{\check67}T_{45}T_{3\check4}T_{\check65}T_{3\hat6}T_{47}\bar T_{1\check4}T_{24}T_{\check 36}T_{\check5\hat8}T_{\check1\check7}T_{\check4\check1}
T_{3\check5}T_{3\hat4}
\underline{T_{5\check1}\bar T_{\check1\check7}}\bar T_{\check 36}
\underline{T_{\check5\hat7}}T_{\hat8\hat7}  
T_{\hat4\check6}T_{2\hat8}  
T_{\check8\check6}
\bar T_{5\hat4}\times \\
\times   
\bar T_{\hat8\hat7}P_{(45\hat4)}=
\zeta^2T_{1\check4}T_{\check67}T_{45}T_{3\check4}T_{\check65}T_{3\hat6}T_{47}\bar T_{1\check4}T_{24}
T_{\check 36}T_{\check5\hat8}\underline{T_{\check1\check7}
T_{\check4\check1}}
T_{3\check5}T_{3\hat4}
\underline{\bar T_{\check1\check7}}T_{5\check1}\underline{P_{(5\check7\hat5)}}
\bar T_{\check 36}
T_{\hat8\hat7}  
T_{\hat4\check6}T_{2\hat8}  
T_{\check8\check6}
\bar T_{5\hat4}\times \\
\times   
\bar T_{\hat8\hat7}P_{(45\hat4)}= 
\zeta^2T_{1\check4}T_{\check67}T_{45}T_{3\check4}T_{\check65}T_{3\hat6}T_{47}\underline{\bar T_{1\check4}T_{24}}T_{\check 36}
T_{\check5\hat8}\underline{T_{\check4\check1}}T_{\check4\check7}
T_{3\check5}T_{3\hat4}
T_{5\check1}
\bar T_{\check 36}T_{\hat4\check6}
T_{\hat85}  
T_{2\hat8}  
T_{\check8\check6}\bar T_{\hat85}
\bar T_{\check7\hat4}\times \\
\times   
P_{(5\check7\hat5)}P_{(45\hat4)}=
\zeta^2T_{1\check4}T_{\check67}T_{45}T_{3\check4}
T_{\check65}T_{3\hat6}\underline{T_{47}}
T_{1\check2}\underline{T_{24}}T_{\check 36}
T_{\check5\hat8}\underline{T_{\check4\check7}}
T_{3\check5}T_{3\hat4}
T_{5\check1}
\bar T_{\check 36}T_{\hat4\check6}
T_{\hat85}  
T_{2\hat8}  
T_{\check8\check6}\bar T_{\hat85}
\bar T_{\check7\hat4}\times \\
\times   
P_{(5\check7\hat5)}P_{(45\hat4)}=
\zeta^3 T_{1\check4}\underline{T_{\check67}T_{45}T_{3\check4}
T_{\check65}T_{3\hat6}}\;
\underline{T_{1\check2}}T_{24}T_{27}P_{(47\hat4)}T_{\check36}
T_{\check5\hat8}
T_{3\check5}T_{3\hat4}
T_{5\check1}
\bar T_{\check 36}T_{\hat4\check6}
T_{\hat85}  
T_{2\hat8}  
T_{\check8\check6}\bar T_{\hat85}
\bar T_{\check7\hat4}\times \\
\times   
P_{(5\check7\hat5)}P_{(45\hat4)}=
\zeta^3 T_{1\check4}T_{1\check2} 
\underline{T_{1\check6}T_{17}
\bar T_{17}\bar T_{1\check6}}
T_{\check67}\underline{T_{24}
\bar T_{24}}T_{45}T_{3\check4}
\underline{T_{\check65}T_{3\hat6}} \; 
\underline{T_{24}}T_{27}P_{(47\hat4)}T_{\check36}
T_{\check5\hat8}
T_{3\check5}T_{3\hat4}
T_{5\check1}
\bar T_{\check 36}T_{\hat4\check6}
T_{\hat85}  
T_{2\hat8}  \times \\
\times T_{\check8\check6}\bar T_{\hat85}
\bar T_{\check7\hat4}   
P_{(5\check7\hat5)}P_{(45\hat4)}=\zeta^3 
\bar\fun_{a_2}
\underline{\bar T_{17}\bar T_{\check67}
\bar T_{1\check6}} T_{24}\underline{\bar T_{24}T_{45}T_{24}}
\; \underline{\bar T_{24}T_{3\check4}T_{24}}
T_{\check65}\underline{T_{3\hat6}
T_{27}}\; \underline{P_{(47\hat4)}T_{\check36}}
T_{\check5\hat8}
T_{3\check5}T_{3\hat4}
T_{5\check1}
\bar T_{\check 36}T_{\hat4\check6}\times \\
\times T_{\hat85}  
T_{2\hat8}  
T_{\check8\check6}\bar T_{\hat85}
\bar T_{\check7\hat4}   
P_{(5\check7\hat5)}P_{(45\hat4)}=\zeta^3 
\bar\fun_{a_2}
T_{\check67}
\bar T_{1\check6} T_{24}T_{25}\underline{T_{45}T_{2\check3}}T_{4\check3}
\underline{T_{\check65}
T_{27}}\underline{T_{3\hat6}T_{\check36}} 
P_{(47\hat4)}
T_{\check5\hat8}
T_{3\check5}T_{3\hat4}
T_{5\check1}
\bar T_{\check 36}T_{\hat4\check6}\times \\
\times T_{\hat85}  
T_{2\hat8}  
T_{\check8\check6}\bar T_{\hat85}
\bar T_{\check7\hat4}   
P_{(5\check7\hat5)}P_{(45\hat4)}=\zeta^4
\bar\fun_{a_2}
\underline{T_{24}T_{25}T_{2\check3}T_{26}} \;  
\underline{\bar T_{26}T_{\check67}}
\bar T_{1\check6}T_{45}T_{4\check3}\underline{T_{27}}
T_{\check65}
P_{(\check6\check36)} 
P_{(47\hat4)}
T_{\check5\hat8}
T_{3\check5}T_{3\hat4}
T_{5\check1}
\bar T_{\check 36}T_{\hat4\check6}\times \\
\times T_{\hat85}  
T_{2\hat8}  
T_{\check8\check6}\bar T_{\hat85}
\bar T_{\check7\hat4}   
P_{(5\check7\hat5)}P_{(45\hat4)}=\zeta^4
\bar\fun_{a_2}\bar \fun_{a_1}
T_{\check67}
\bar T_{26}\bar T_{1\check6}
T_{45}T_{4\check3}
\underline{T_{\check65}
T_{\check5\hat8}
T_{\hat6\check5}} \; 
\underline{T_{6\hat7}}T_{5\check1}\underline{\bar T_{63}T_{\hat7\check3}}T_{\hat85} T_{2\hat8}  
T_{\check8\check3}\bar T_{\hat85}
\bar T_{4\hat7}\times \\
\times    
P_{(\check6\check36)} 
P_{(47\hat4)}
P_{(5\check7\hat5)}P_{(45\hat4)}= 
\zeta^6
\bar\fun_{a_2}\bar \fun_{a_1}
T_{\check67}
\bar T_{26}
\underline{\bar T_{1\check6}}
T_{45}T_{4\check3}
T_{\check5\hat8}
\underline{T_{\check6\hat8}
T_{6\check1}}
\bar T_{\hat53}
T_{\check5\hat7}
T_{\hat86} 
T_{2\hat8}  
T_{\check8\check7}
\bar T_{\hat86}
\bar T_{4\bar3}
\times \\
\times P_{(\check656)}
P_{(\hat3\hat7\check3)}   
P_{(\check6\check36)} 
P_{(47\hat4)}
P_{(5\check7\hat5)}P_{(45\hat4)}=
\zeta^6
\bar\fun_{a_2}\bar \fun_{a_1}
T_{\check67}
\bar T_{26}
T_{45}T_{4\check3}
T_{\check5\hat8}
T_{1\hat8}
\underline{T_{\check6\hat8}}
\bar T_{\hat53}
T_{\check5\hat7}
\underline{T_{\hat86}} 
T_{2\hat8}  
T_{\check8\check7}
\bar T_{\hat86}
\bar T_{4\bar3}
\times \\
\times P_{(\check656)}
P_{(\hat3\hat7\check3)}   
P_{(\check6\check36)} 
P_{(47\hat4)}
P_{(5\check7\hat5)}P_{(45\hat4)}=
\zeta^7
\bar\fun_{a_2}\bar \fun_{a_1}
T_{\check67} 
\bar T_{26}
T_{45}T_{4\check3}
T_{\check5\hat8}
T_{1\hat8}
\bar T_{\hat53}
T_{\check5\hat7} P_{(\check6\hat86)}
T_{\hat8} 
T_{\check8\check7}
\bar T_{\hat86}
\bar T_{4\bar3}
\times \\
\times 
P_{(\check656)}
P_{(\hat3\hat7\check3)}   
P_{(\check6\check36)} 
P_{(47\hat4)}
P_{(5\check7\hat5)}P_{(45\hat4)}=
\zeta^7
\bar\fun_{a_2}\bar \fun_{a_1}
\underline{T_{\check67}} \; 
\underline{\bar T_{26}}
T_{45}T_{4\check3}
T_{\check5\hat8}
T_{1\hat8}
\bar T_{\hat53}
\underline{T_{\check5\hat7}} \; 
\underline{T_{26}} \; \underline{T_{\hat6\check7}}
\bar T_{6\check8}
\bar T_{4\bar3}
\times \\
\times P_{(\check6\hat86)}P_{(\check656)}
P_{(\hat3\hat7\check3)}   
P_{(\check6\check36)} 
P_{(47\hat4)}
P_{(5\check7\hat5)}P_{(45\hat4)}=
\zeta^8
\bar\fun_{a_2}\bar \fun_{a_1}
T_{45}T_{4\check3}
T_{\check5\hat8}
T_{1\hat8}
\bar T_{\hat53}
T_{\check5\hat7}T_{\check6\hat5}P_{(\check676)}
\bar T_{6\check8}
\bar T_{4\bar3}
\times \\
\times P_{(\check6\hat86)}P_{(\check656)}
P_{(\hat3\hat7\check3)}   
P_{(\check6\check36)} 
P_{(47\hat4)}
P_{(5\check7\hat5)}P_{(45\hat4)}=
\zeta^8
\bar\fun_{a_2}\bar \fun_{a_1}\bar \fun_{a_0}
\bar T_{\check8\check7}\bar T_{\check8\check1}\bar T_{\check8\check 4}
\underline{\bar T_{\check8\hat5}
T_{45}}T_{4\check3}
\underline{T_{\check5\hat8}}
T_{1\hat8}
\bar T_{\hat53}
T_{\check5\hat7}T_{\check6\hat5}
\bar T_{\hat7\check8}
\bar T_{4\bar3}
\times \\
\times P_{(\check676)}P_{(\check6\hat86)}P_{(\check656)}
P_{(\hat3\hat7\check3)}   
P_{(\check6\check36)} 
P_{(47\hat4)}
P_{(5\check7\hat5)}P_{(45\hat4)}=
\zeta^8
\bar\fun_{a_2}\bar \fun_{a_1}\bar \fun_{a_0}
\bar T_{\check8\check7}\underline{\bar T_{\check8\check1}} \; \underline{\bar T_{\check8\check 4}} \; 
\underline{T_{4\hat8}} T_{45}
T_{4\check3}
\underline{T_{1\hat8}}
\bar T_{\hat53}
T_{\check5\hat7}T_{\check6\hat5}
\bar T_{\hat7\check8}
\bar T_{4\bar3}
\times \\
\times P_{(\check676)}P_{(\check6\hat86)}P_{(\check656)}
P_{(\hat3\hat7\check3)}   
P_{(\check6\check36)} 
P_{(47\hat4)}
P_{(5\check7\hat5)}P_{(45\hat4)}=
\zeta^8
\bar\fun_{a_2}\bar \fun_{a_1}\bar \fun_{a_0}
\underline{\bar T_{\check8\check7}} T_{45}
T_{4\check3}
\bar T_{\hat53}
\underline{T_{\check5\hat7}}T_{\check6\hat5}
\underline{\bar T_{\hat7\check8}}
\bar T_{4\bar3}
\times \\
\times P_{(\check676)}P_{(\check6\hat86)}P_{(\check656)}
P_{(\hat3\hat7\check3)}   
P_{(\check6\check36)} 
P_{(47\hat4)}
P_{(5\check7\hat5)}P_{(45\hat4)}
\end{multline*}
In the previous lines we used both the Pentagon relation 
coupled with the symmetry property several times and the commutativity 
relations corresponding to the underlined fragments. 
Sometimes several simplifications are recorded  
in the same line, as in the first equality above where the underlined 
factors $\bar{T}_{24}$ and $T_{24}$ commute with $T_{\hat8\hat7}T_{\check5\hat8}$ 
and therefore cancel each other, so that along with the 
first underlined factor we obtain  a subproduct 
$\bar{T}_{\check5\hat8}T_{\hat8\hat7}T_{\check5\hat8}$ and the Pentagon 
relation can be applied.

\vspace{0.2cm}\noindent 
Use now the identity: 
\[ \bar T_{\check8\check7}T_{\check5\hat7}\bar T_{\hat7\check8}=
T_{\check5\hat8}T_{\check5\hat7}\bar T_{\check8\check7}\bar T_{\hat7\check8}=
\zeta^{-1}T_{\check5\hat8}T_{\check5\hat7}P_{(\check7\check8\hat7)}\]
and introduce above to find that: 
\begin{multline*}
\bar\fun_{a_{12}}\bar\fun_{a_{23}}\bar\fun_{a_{13}}=
\zeta^7
\bar\fun_{a_2}\bar \fun_{a_1}\bar \fun_{a_0}
T_{45}
T_{4\check3}
\bar T_{\hat53}
T_{\check5\hat8}
T_{\check5\hat7}
T_{\check6\hat5}
\bar T_{4\bar3}
\times \\
\times P_{(\check7\check8\hat7)}
P_{(\check676)}P_{(\check6\hat86)}P_{(\check656)}
P_{(\hat3\hat7\check3)}   
P_{(\check6\check36)} 
P_{(47\hat4)}
P_{(5\check7\hat5)}P_{(45\hat4)}=\zeta^7
\bar\fun_{a_2}\bar \fun_{a_1}\bar \fun_{a_0}\bar \fun_{a_3}
\bar T_{3\hat6}\bar T_{3\hat7}\bar T_{3\hat8}\underline{\bar T_{3\check5}
T_{45}T_{4\check3}}
\bar T_{\hat53}
T_{\check5\hat8}
T_{\check5\hat7}
T_{\check6\hat5}
\bar T_{4\bar3}
\times \\
\times P_{(\check7\check8\hat7)}
P_{(\check676)}P_{(\check6\hat86)}P_{(\check656)}
P_{(\hat3\hat7\check3)}   
P_{(\check6\check36)} 
P_{(47\hat4)}
P_{(5\check7\hat5)}P_{(45\hat4)}=\zeta^7
\bar\fun_{a_2}\bar \fun_{a_1}\bar \fun_{a_0}\bar \fun_{a_3}
\bar T_{3\hat6}\bar T_{3\hat7}\bar T_{3\hat8}
T_{45}
\underline{\bar T_{3\check5}
\bar T_{\hat53}}
T_{\check5\hat8}
T_{\check5\hat7}
T_{\check6\hat5}
\bar T_{4\bar3}
\times \\
\times P_{(\check7\check8\hat7)}
P_{(\check676)}P_{(\check6\hat86)}P_{(\check656)}
P_{(\hat3\hat7\check3)}   
P_{(\check6\check36)} 
P_{(47\hat4)}
P_{(5\check7\hat5)}P_{(45\hat4)}= 
\zeta^6
\bar\fun_{a_2}\bar \fun_{a_1}\bar \fun_{a_0}\bar \fun_{a_3}
\bar T_{3\hat6}\bar T_{3\hat7}\bar T_{3\hat8}
T_{45}
P_{(\check53\hat5)}
T_{\check5\hat8}
T_{\check5\hat7}
T_{\check6\hat5}
\bar T_{4\bar3}
\times \\
\times P_{(\check7\check8\hat7)}
P_{(\check676)}P_{(\check6\hat86)}P_{(\check656)}
P_{(\hat3\hat7\check3)}   
P_{(\check6\check36)} 
P_{(47\hat4)}
P_{(5\check7\hat5)}P_{(45\hat4)}=
\zeta^6
\bar\fun_{a_2}\bar \fun_{a_1}\bar \fun_{a_0}\bar \fun_{a_3}
\underline{\bar T_{3\hat6}}\; \underline{\bar T_{3\hat7}} \; \underline{\bar T_{3\hat8}} \; \underline{
T_{45}} \; 
\underline{T_{3\hat8}} \; \underline{T_{3\hat7}} \underline{T_{\check6\check3}} \;
\underline{\bar T_{45}}
\times \\
\times P_{(\check53\hat5)}P_{(\check7\check8\hat7)}
P_{(\check676)}P_{(\check6\hat86)}P_{(\check656)}
P_{(\hat3\hat7\check3)}   
P_{(\check6\check36)} 
P_{(47\hat4)}
P_{(5\check7\hat5)}P_{(45\hat4)}=\zeta^6
\bar\fun_{a_2}\bar \fun_{a_1}\bar \fun_{a_0}\bar \fun_{a_3}
\end{multline*}
Thus the lift of the lantern relation is $\zeta^6$. 
Therefore we have to renormalize each right Dehn twist 
by taking $\widetilde{D_{\alpha}}=\zeta^{-6}\fun_{\alpha}$, as claimed. 
\end{proof}

\vspace{0.1cm}\noindent 
The following lemma is a simple consequence of a deep result of Gervais
from (\cite{Ge}):
\begin{lemma}\label{gervais}
Let $g\geq 2$ and $s\geq 0$. Then the 
group $\Gamma^s_{g,r}$ is presented as follows:
\begin{enumerate}
\item Generators are all Dehn twists $D_{a}$ along the  non-separating
simple closed curves $a$ on $\Sigma^s_{g,r}$.
\item Relations:
\begin{enumerate}
\item Braid type 0-relations:
\[ D_aD_b=D_bD_a\]
for each pair of disjoint non-separating simple closed curves $a$ and $b$;
\item   Braid type 1-relations:
\[ D_aD_bD_a=D_bD_aD_b\]
for each pair of non-separating simple closed curves
$a$ and $b$ which intersect transversely in one point;
\item One lantern relation for a $4$-hold sphere embedded in $\Sigma^{s}_{g,r}$ so that all
boundary curves are non-separating;
\item One chain relation for a 2-holed torus embedded in $\Sigma^s_{g,r}$ so that all
boundary curves are non-separating;
\item A puncture relation for each puncture.
\end{enumerate}
\end{enumerate}
\end{lemma}
\begin{proof} 
According to (\cite{Ge}, Theorem B) we have a presentation
of $\Gamma_{g,s+r}$  with the generators above and all but the puncture relations.
Now, the kernel of $\Gamma_{g,s+r}\to \Gamma^s_{g,r}$ is  the free Abelian
group generated  by the Dehn twists  along the boundary curves to be pinched to punctures.
Such a Dehn twist is expressed (using the lantern relation) by the left hand side
of the puncture relation. This proves the claim.
\end{proof}

\vspace{0.2cm}\noindent
{\em Proof of Proposition \ref{pres}}. According to the normalization
coming from the braid relations and the lantern relations the
images of the standard  Dehn twist
generators of the mapping class group are products of $\zeta^6$ and
elements $T_{ij}$, where $i,j$ are the labels of the triangles (possibly
with $\:\hat{}\:$ or $\:\check{}\:$). 
Thus the projective factors that appear belong to the
subgroup $A$ generated by $\zeta^6$. The only non-trivial lift of a
relation from Lemma \ref{gervais} is the chain relation which lifts to
$\zeta^{-72}$. Set $z$ for the element $\zeta^{-6}$ of 
$\widetilde{\Gamma_{g,r}^s}$. Then the presentation of the central
extension $\widetilde{\Gamma^s_{g,r}}$ is given by the claimed relations.

\subsection{Cohomological consequences}
Recall from (\cite{KS}, Corollary 4.4) 
that the 2-cohomology classes $\chi$ and $e_i$
are defined for any $g\geq 3, s,r\geq 0$ and they span a free Abelian 
subgroup $\Z^{s+1}\subset H^2(\Gamma_{g,r}^s)$.  
This inclusion is actually an isomorphism when $g\geq 4$.

\vspace{0.2cm}\noindent
We will denote by 
$\widehat{\Gamma_{g,r}^s}$ the group defined by the presentation 
given in Proposition \ref{pres}, for all values of $s,g,r$. 
Thus, according to Proposition \ref{pres}  the extension 
$\widehat{\Gamma_{g,r}^s}$ is isomorphic to 
$\widetilde{\Gamma_{g,r}^s}$ if $s\geq 4$ and $g\geq 2$. 

\begin{lemma}
If $g\geq 2$, then we have  
$c_{\widehat{\Gamma_{g,r}}}=12\chi\in H^2(\Gamma_{g,r};A)$.
\end{lemma}
\begin{proof}
Consider first the case where $\zeta$ is not a root of unity, so 
that the group $A$ is isomorphic to $\Z$. 
Gervais proved in (\cite{Ge}, Theorem 3.6) that  
$\widehat{\Gamma_{g,r}}$ (namely,  
where $s=0$) is isomorphic to the so-called $p_1$-central extension of 
$\Gamma_{g,r}$. Further in \cite{Ge,MR} the authors identified 
the class of the $p_1$-central extension of $\Gamma_{g,r}$ 
to the class $12\chi$ and thus $c_{\widehat{\Gamma_{g,r}}}=12\chi$.

\vspace{0.2cm}
\noindent
Here is a more direct argument. 
Set $\Gamma_{g,r}(1)$  for the subgroup of
$\widehat{\Gamma_{g,r}}$ generated by the lifts
$\widetilde{D_{a}}$ of the Dehn twists and 
the central element $u=z^{12}$.  Then
$\Gamma_{g,r}(1)$ is the universal central extension considered by 
Harer (see \cite{Ge,Harer}) and thus 
$c_{\Gamma_{g,r}(1)}$ is the generator
$\chi$ of $H^2(\Gamma_{g,r})\cong \Z$.

\vspace{0.2cm}
\noindent
The cohomology class $c_{\Gamma_{g,r}(1)}$ is represented 
by some explicit 2-cocycle $C_{\Gamma_{g,r}(1)}:\Gamma_{g,r}\times \Gamma_{g,r}\to \Z$ which arises as follows. Let 
$S:\Gamma_{g,r}\to \Gamma_{g,r}(1)$ be a set-wise section. 
Let also $i:\ker(\Gamma_{g,r}(1)\to \Gamma_{g,r})\to \Z$ 
be the group isomorphism defined by $i(u)=1$. 
It is well-known that the  2-cocycle
\[ C_{\Gamma_{g,r}(1)}(x,y)= i(S(xy)S(x)^{-1}S(y)^{-1}) \in\Z \]
represents the cohomology class $c_{\Gamma_{g,r}(1)}$. 

\vspace{0.2cm}
\noindent 
Let us construct now a 2-cocycle representing 
the extension $\widehat{\Gamma_{g,r}}$. Consider the set-wise section 
$\iota\circ S:\Gamma_{g,r}\to \widehat{\Gamma_{g,r}}$, where 
$\iota:\Gamma_{g,r}(1)\to \widehat{\Gamma_{g,r}}$ is the  obvious 
inclusion. Let also $j: \ker(\widehat{\Gamma_{g,r}}\to 
\Gamma_{g,r})\to \Z$ be the isomorphism given by $j(z)=1$. 
Then 
\[ C_{\widehat{\Gamma_{g,r}}}(x,y)= j((\iota\circ S)(xy)
(\iota\circ 
S)(x)^{-1}(\iota\circ S)(y)^{-1}) =j(\iota(S(xy)S(x)^{-1}S(y)^{-1})) \in\Z \]
is a 2-cocycle representing $c_{\widehat{\Gamma_{g,r}}}$.
Since $j(\iota(u))=j(z^{12})=12 i(u)$ and 
$S(xy)S(x)^{-1}S(y)^{-1}$ belongs to the cyclic subgroup of 
$\Gamma_{g,r}(1)$ generated by $u$,  it follows that 
\[ C_{\widehat{\Gamma_{g,r}}}(x,y)= 12 C_{\Gamma_{g,r}(1)}\]
and thus $c_{\widehat{\Gamma_{g,r}}}=12\chi$, where 
$\chi$ is one fourth of the  Meyer signature class, which is a
generator of $H^2(\Gamma_{g,1})\subset H^2(\Gamma^1_g)$.

\vspace{0.2cm}
\noindent When $\zeta$ is a root of unity of order $N$ then 
the class of the extension $\widehat{\Gamma_{g,r}}$ is the image 
of $12\chi$ in $H^2(\Gamma_{g,r};\Z/N\Z)$ by the reduction mod $N$.  
\end{proof}

\vspace{0.2cm}
\noindent
The next step is to prove a similar statement when the number
$s$ of punctures is non-zero.

\begin{definition}\label{extens}
For $(m_1,m_2,\ldots,m_s)\in\Z^s$ let 
$\Gamma^s_{g,r}(m_1,m_2,\ldots,m_s)$ be the central extension 
of $\Gamma_{g,r}^s$ by $A$ having  the following presentation: 
\begin{enumerate}
\item Generators are the $\widetilde{D_{\alpha}}$, where 
$D_{\alpha}$ are Dehn twist generators of 
$\Gamma^s_{g,r}$ and the central element $z$ of the same order as $\zeta^{-6}$;
\item Relations are as follows. For each puncture $p_i$  the 
lift of the corresponding puncture 
relation reads: 
\[ \widetilde{D_{a_1(i)}}^{-1}\widetilde{D_{a_2(i)}}^{-1}\widetilde{D_{a_3(i)}}^{-1}
\widetilde{D_{a_{12}(i)}}\widetilde{D_{a_{13}(i)}}\widetilde{D_{a_{23}(i)}}=z^{m_i}\]
where $\widetilde{D_a}$ are lifts of Dehn twists.
Furthermore the chain and lantern relations have trivial lifts. 
\end{enumerate}
\end{definition}

\begin{proposition}\label{euler}
Suppose that $g\geq 0$. 
Then $c_{\Gamma^s_{g,r}(m_1,\ldots,m_s)}\in A^{n+1}\subset H^2(\Gamma^s_{g,r};A)$
is the vector $m_1e_1+m_2e_2+\cdots+m_se_s$, where
$e_i$ is the Euler class of the $i$-th puncture.
\end{proposition}
\begin{proof}
This is folklore. Consider first that $\zeta$ is not a root of unity. 
Let $\Sigma^{s-1}_{g,r+1;i}$ denote the subsurface 
of $\Sigma_{g,r}^s$ obtained by removing a one-punctured disk 
centered at the puncture $p_i$ and thus creating a new boundary 
component $b_i$.  
We have then a central extension 
\[ \Z\to \Gamma^{s-1}_{g,r+1;i}\to \Gamma^s_{g,r}\to 1\]
induced by  the  inclusion map  
$\Sigma^{s-1}_{g,r+1;i}\hookrightarrow \Sigma^s_{g,r}$. 
It is well-known that its cohomology class is 
$c_{\Gamma^{s-1}_{g,r+1;i}}=e_i$.

\begin{lemma}\label{lifts}
The extension $\Gamma^{s-1}_{g,r+1;i}$ is isomorphic to 
$\Gamma_{g,r}^s(0,\ldots,1,0\ldots,0)$, where 1 is on the $i$-th position. 
\end{lemma}
\begin{proof}
There is a natural set-wise section 
$S_i:\Gamma^s_{g,r}\to \Gamma^{s-1}_{g,r+1;i}$, given by 
$S_i(D_{\alpha})=D_{\alpha}$, for any Dehn twist $D_{\alpha}$. 
In order to make sense, we might suppose that a simple 
closed curve $\alpha$ disjoint from the 
puncture $p_i$ is actually disjoint from $b_i$ 
so that it lies within $\Sigma^{s-1}_{g,r+1;i}$. 

\vspace{0.1cm}
\noindent 
Braid, chain and lantern relations 
are then lifted trivially. A puncture relation at $p_j$ is lifted trivially 
if $j\neq i$. Consider next a puncture relation at 
$p_i$ in  $\Sigma_{g,r}^s$, which 
is supported on some subsurface $\Sigma^{1}_{0,3}$. The three boundary 
curves of $\Sigma^{1}_{0,3}$ lie within $\Sigma^{s-1}_{g,r+1;i}$  
and together with $b_i$ bound  a 4-holed sphere 
in $\Sigma^{s-1}_{g,r+1;i}$. The lantern relation associated to 
this 4-holed sphere on  $\Sigma^{s-1}_{g,r+1;i}$ is then the  
lift of the puncture relation at $p_i$.  The 
Dehn twist along $b_i$ is the generator $z$ of the central factor 
$\ker(\Gamma^{s-1}_{g,r+1;i}\to \Gamma^s_{g,r})$. 
Thus the lift of a puncture relation at $p_i$ is the factor $z$. 
\end{proof}

\begin{lemma}
Let $L_{\mathbf m}:\Z^s\to\Z$ denote the linear map 
$L_{\mathbf m}(n_1,\ldots,n_s)=\sum_{i=1}^s m_in_i$, where 
${\mathbf m}=(m_1,\ldots,m_s)$. Consider the central extension 
\[ 1\to \Z^s\to \Gamma_{g,r+s}\to \Gamma_{g,r}^s\to 1\] 
Then the map $L_{\mathbf m}$ induces a quotient of $\Gamma_{g,r+s}$, which 
is a central extension $\Gamma_{g,r}^s({\bf m})$ of $\Gamma_{g,r}^s$  
by $\Z$ which is isomorphic to $\Gamma_{g,r}^s(m_1,m_2,\ldots,m_s)$ and 
gives rise to the following commutative diagram: 
\[\begin{array}{cclcclccc} 
1 &\to& \Z^s&\to& \Gamma_{g,r+s}&\to &\Gamma_{g,r}^s&\to& 1\\
 & & \downarrow L_{\bf a} & & \downarrow  \pi  & & \downarrow {\bf 1}  & \\
1 &\to& \Z&\to& \Gamma_{g,r}^s({\bf m}) &\to &\Gamma_{g,r}^s&\to& 1\\
\end{array}
\]
\end{lemma}
\begin{proof}
The class of the central extension $c_{\Gamma_{g,r+s}}$ belongs to 
$H^2(\Gamma_{g,r}^s;\Z^s)= 
\oplus_{s}H^2(\Gamma_{g,r}^s,\Z)$. By functoriality we derive that 
$c_{\Gamma_{g,r+s}}=(e_1,e_2,\ldots,e_s)\in H^2(\Gamma_{g,r}^s;\Z^s)$. 
Then the class $c_{\Gamma_{g,r}^s({\bf a})}$ 
is the image of $c_{\Gamma_{g,r+s}}$ into $H^{2}(\Gamma_{g,r}^s)$ 
by the homomorphism of coefficients rings 
$L_{\bf m}:\Z^s\to\Z$. There is an obvious set-wise section $S$ defined  
in the same way as the $S_i$ from above. Then $c_{\Gamma_{g,r+s}}$ 
is the class of  the 2-cocycle $L_{\bf m}C$, where $C$ is the 2-cocycle 
associated to $S$ and so 
\[ L_{\bf m}C(x,y)=\pi(S(x)^{-1}S(y)^{-1}S(xy))=
L_{\bf m}((S_i(x)^{-1}S_i(y)^{-1}S_i(xy)_{i=1,s})=
\sum_{i=1}^sm_iC_i(x,y)\]
where $C_i$ is the 2-cocycle associated to $S_i$. 
Since the class of $C_i$  is $e_i$ it follows that the 
class of $L_{\bf m}C$ is $\sum_{i=1}^s m_iei$. 

\vspace{0.2cm}\noindent
On the other hand the lifts of relations in $\Gamma_{g,r}^s({\bf m})$ 
are the same as in $\Gamma_{g,r}^s(m_1,\ldots,m_s)$ and thus they 
are isomorphic. In fact the lifts of braid, chain and lantern relations 
to $\Gamma_{g,r+s}$ are trivial. The lift of a puncture relation 
at $p_i$ is the $i$-th generator of the central factor $\Z^s$, 
according to Lemma \ref{lifts}. Therefore its image into 
$\Gamma_{g,r}^s({\bf m})$ is $z^{m_i}$, namely the lift of the puncture 
relation in $\Gamma_{g,r}^s(m_1,\ldots,m_s)$. 
\end{proof}

\vspace{0.1cm}
\noindent When $\zeta$ is a root of unity the extensions 
by $\Z$ above are replaced by extensions by $\Z/N\Z$ and all arguments go 
through without essential modifications.

\vspace{0.1cm}
\noindent This proves the Proposition.
\end{proof}

\vspace{0.2cm}\noindent
{\em Proof of  the Theorem.}
Assume first that $A$ is cyclic infinite. Consider the  
operation $\otimes$ (which is a push-out, or a fibered product) 
on central extensions defined 
as follows. If $f_i:G_i\to G$ are the projections 
homomorphisms of the central extensions $G_i$ of $G$ by $\Z$ then 
$G_1\otimes G_2$ is the extension $f_1^*G_2$ (or equivalently $f_2^*G_1$) 
of $G$ by $\Z^2$. The class $c_{G_1\otimes G_2}\in H^2(G,\Z^2)$ 
is the direct sum of the classes $c_{G_i}\in  H^2(G,\Z)$ under the 
identification of $H^2(G,\Z^2)$ with the sum of two copies of $H^2(G,\Z)$.

\vspace{0.2cm}\noindent 
Let $f$ denote the surjective homomorphism
$f:\Gamma_{g,r}^s\to \Gamma_{g,r}$. 
Consider then the central extension 
\[ 1\to \Z^2\to f^*(\widehat{\Gamma_{g,r}})\otimes 
\Gamma_{g,r}^s(1,1,\ldots,1)\to \Gamma_{g,r}^s\to 1\]
Using the map $L:\Z^2\to \Z$ given by 
$L(x,y)=x+y$ we find a quotient of $f^*(\widehat{\Gamma_{g,r}})\otimes 
\Gamma_{g,r}^s(1,1,\ldots,1)$, which is a central extension by $\Z$ isomorphic  
to $\widehat{\Gamma_{g,r}^s}$. 
In fact, there is a commutative diagram:  
\[\begin{array}{cclcclccc} 
1 &\to& \Z^2&\to& f^*(\widehat{\Gamma_{g,r}})\otimes 
\Gamma_{g,r}^s(1,1,\ldots,1)&\to &\Gamma_{g,r}^s&\to& 1\\
 & & \downarrow L & & \downarrow  \pi  & & \downarrow {\bf 1}  & \\
1 &\to& \Z&\to& \widehat{\Gamma_{g,r}^s} &\to &\Gamma_{g,r}^s&\to& 1\\
\end{array}
\]
The central extension from the lower row is isomorphic to 
$\widehat{\Gamma_{g,r}^s}$ because the lifts of relations are the same. 
Braid and lantern relations lift trivially. 
Chain relations lift to $z^{12}$ in  $f^*(\widehat{\Gamma_{g,r}})$ and 
trivially to $\Gamma_{g,r}^s(1,1,\ldots,1)$ and thus the image of the lift 
by $L$ (or $\pi$) is $z^{12}$. Puncture relations at $p_i$ lift trivially 
to $f^*(\widehat{\Gamma_{g,r}})$ and to $z$ in the factor 
$\Gamma_{g,r}^s(1,1,\ldots,1)$, so that its image by $L$ (or $\pi$) is 
$z$. 
As a consequence of this description the class 
$c_{\widehat{\Gamma_{g,r}^s}}$ is the image by $L$ of the 
class of $f^*(\widehat{\Gamma_{g,r}})\otimes 
\Gamma_{g,r}^s(1,1,\ldots,1)$, namely 
$c_{f^*(\widehat{\Gamma_{g,r}})}+c_{\Gamma_{g,r}^s(1,1,\ldots,1)}$.

\vspace{0.2cm}\noindent 
On the other hand, by functoriality, 
the class $c_{f^*(\widehat{\Gamma_{g,r}})}$ is 
$f^*(12\chi)=12\chi\in H^{2}(\Gamma_{g,r}^s)$, because 
the map $f^*$ is the standard embedding 
of $H^2(\Gamma_{g,r})=\Z\chi$ into 
$H^2(\Gamma_{g,r}^s)$. Proposition \ref{euler} proves the Theorem 
for $g\geq 3$. 

\vspace{0.2cm}\noindent 
When $g=2$ one does not know the group $H^2(\Gamma_{2,r}^s)$, but for $s=0$ and 
$r\leq 1$. Nevertheless, the classes $\chi$ and $e_j$ are still defined. 
It suffices to prove that: 

\begin{lemma}
The subgroup of $H^2(\Gamma_{2,r}^s)$ 
generated by $\chi$ and $e_1,\ldots,e_s$ is isomorphic to 
$\Z/10\Z\oplus \Z^s$.  
\end{lemma}
\begin{proof}
By the universal coefficients theorem we have 
\[1\to H_1(\Gamma_{2,r}^s)\to H^2(\Gamma_{2,r}^s)\to {\rm Hom}(H_2(\Gamma_{2,r}^s), \Z)\to 1\]
From (\cite{KS}, Proposition 1.6) we have $H_1(\Gamma_{2,r}^s)=\Z/10\Z$. 
The Meyer class $\chi$ in genus 2 is one half of the class 
of Meyer's cocycle from \cite{Me} and it generates the image 
of $H_1(\Gamma_{2,r}^s)$ into  $H^2(\Gamma_{2,r}^s)$. 

\vspace{0.2cm}\noindent 
Consider next the extensions $\Gamma_{2,r}^s({\bf m})$ for 
integral vectors ${\bf m}$. According to the previous description 
lifts of puncture relations are of the form $z^{m_i}$. Suppose that 
there exists an isomorphism between the extensions 
$\Gamma_{2,r}^s({\bf m})$ and  $\Gamma_{2,r}^s({\bf u})$. 
Such an isomorphism of extensions 
should send $\widetilde{D_{\alpha}}$ into $z^{n(\alpha)}
\widetilde{D_{\alpha}}$, because it  has to induce 
the identity on $\Gamma_{2,r}^s$. Since lifts of braid relations 
are trivial in both extension groups it follows that 
$n(\alpha)=n$ does not depend on the non-separating curve $\alpha$.  
But puncture relations are homogeneous, and so they  
do not depend on $n$. This shows that 
${\bf m}={\bf u}$. In particular the classes $e_i$ span a free 
$\Z$-submodule of $H^2(\Gamma_{2,r}^s)$.   

\vspace{0.2cm}\noindent 
Since the class $\chi$ is of order 10 and both subgroups 
$\Z/10\Z$ (generated by $\chi$) and $\Z^s$ (generated by $e_1,\ldots,e_s$)
inject into  $H^2(\Gamma_{2,r}^s)$, the claim follows. 
\end{proof}

\vspace{0.2cm}\noindent 
Then the arguments used above for $g\geq 3$ work as well for $g=2$ and 
the Theorem follows. When $\zeta$ is a root of unity 
the associated cohomology class is the reduction mod $N$ 
of the corresponding integral cohomology class.

\vspace{0.2cm}\noindent 
{\em Proof of Corollary \ref{cor1}}. 
Consider the extension $\widehat{\Gamma_{g,r+s}}$ of class 
$12\chi$. The Corollary claims that there is an exact sequence: 
\[ 1\to A^{s-1}\to \widehat{\Gamma_{g,r+s}}\to \widehat{\Gamma_{g,r}^s}\to 1\] 
This can be verified by using the explicit presentations of the two groups 
involved. The kernel is generated by the products of two opposite Dehn 
twists on the $s$ blown up boundary components.

\vspace{0.2cm}\noindent 
{\em Proof of Corollary \ref{cor2}}.
It suffices to understand the 
map $H^2(\Gamma_{g,r}^s;A)\to H^2(\Gamma_{g,r}^s,\C^*)$ 
induced by $z\to \zeta^{-6}$. This map is injective,  when $g\geq 3$.

\vspace{0.1cm}\noindent 
The Universal Coefficients Theorem states that, for any Abelian group $W$, 
the following exact sequence is exact: 
\[1\to{\rm Ext}(H_0(\Gamma_{g,r}^s), W)\to H^1(\Gamma_{g,r}^s; W)\to 
{\rm Hom}(H_1(\Gamma_{g,r}^s), W)\to 1\]
Now ${\rm Ext}(\Z,W)=0$, for any  Abelian group $W$.  This implies that 
$H^1(\Gamma_{g,r}^s; \C^*)=H^1(\Gamma_{g,r}^s; \C^*/A)=0$, 
if $g\geq 3$. 
From the Bockstein exact sequence
\[H^1(\Gamma_{g,r}^s;\C^*)\to 
H^1(\Gamma_{g,r}^s; \C^*/A)\stackrel{\beta}{\to} H^2(\Gamma_{g,r}^s;A)\stackrel{\nu}{\to} H^2(\Gamma_{g,r}^s;\C^*)\]
we derive the claim. 

\vspace{0.1cm}\noindent
When  $g=2$ the Universal Coefficient Theorem 
shows, as above,  
that $H^1(\Gamma_{2,r}^s; \C^*)={\rm Hom}(H_1(\Gamma_{2,r}^s), \C^*)$ 
and  $H^1(\Gamma_{2,r}^s; \C^*/A)={\rm Hom}(H_1(\Gamma_{2,r}^s), \C^*/A)$.  
Thus $H^1(\Gamma_{2,r}^s; \C^*)={\rm Hom}(\Z/10\Z,\C^*)=U_{10}$, where 
$U_{10}$ is the subgroup of roots of unity of order $10$. The last 
isomorphism sends a homomorphism into its value on the generator $1$. 
Next $H^1(\Gamma_{2,r}^s; \C^*/A)={\rm Hom}(\Z/10\Z, \C^*/A)=U_{10}\times A/10A$. 
To explain the last isomorphism, each element 
$f\in H^1(\Gamma_{2,r}^s; \C^*/A)$  is determined by its value 
$f(1)=As$, for some $s\in\C^*$. Here $s^{10}=a^n\in A$, where $a$ is the generator of $A$.   Fix some $10$-th root $a^{1/10}\in \C^*$ of the generator of 
$A$. Then the isomorphism above associates to 
$f$ the element $(sa^{-n/10}, s^{10})\in 
U_{10}\times A/10A$, which is well-defined and independent of the choice of 
the representative $s$ in its $A$-coset. 
In particular the map $H^1(\Gamma_{2,r}^s, \C^*)\to H^1(\Gamma_{2,r}^s, \C^*/A)$ 
sends $U_{10}$ onto the factor $U_{10}$ of the second group.  

\vspace{0.1cm}\noindent  
Let $\widehat{f}$ be a lift of $f$ to 
$\widehat{f}:\Z/10\Z=H_1(\Gamma_{g,r}^s)\to\C^*$, for instance 
$\widehat{f}(k)=s^k$, where $k\in \Z/10\Z$. Then 
$F(k_1,k_2)=\widehat{f}(k_1)\widehat{f}(k_2)\widehat{f}(k_1k_2)^{-1}\in A$ 
is a 2-cocycle on $H_1(\Gamma_{2,r}^s)$ with values in $A$. 
The pull-back in $H^2(\Gamma_{2,r}^s,A)$ of the class of 
$F$ by the map $\Gamma_{2,r}^s\to H_1(\Gamma_{2,r}^s)$ is 
the element $\beta(f)$. It is well-known that 
$H^2(\Z/10\Z,A)=A/10A$ is generated by the Euler class. 
Specifically, the cohomology class of the 2-cocycle $F$ in 
$H^2(\Z/10\Z,A)$ is the element $s^{10}\in A/10A$, under the previous 
isomorphism. 

\vspace{0.1cm}\noindent
The Universal Coefficients Theorem shows that 
\[ 1\to {\rm Ext}(H_1(\Gamma_{2,r}^s), A)\to  
H^2(\Gamma_{2,r}^s;A) \to {\rm Hom}(H_2(\Gamma_{2,r}^s), A)\to 1\]
Further ${\rm Ext}(H_1(\Gamma_{g,r}^s), A)=A/10A$ is generated by 
the class $\chi$ (as an $A$-valued cohomology class).  
Using the definition of {\rm Ext} one identifies the class $\chi$ with 
the generator of $H^2(\Z/10\Z; A)$. This implies that the image of $\beta$ is 
the subgroup generated by $\chi$ within $H^2(\Gamma_{2,r}^s; A)$. 
Then Corollary \ref{cor2} follows.

{
\small

\bibliographystyle{plain}

}

\end{document}